\documentclass[11pt]{amsart}

\usepackage{fullpage}
\usepackage{url}
\usepackage{amssymb}
\usepackage{cite}

\renewcommand{\a}{\alpha}

\newcommand{\g}{\gamma}
\renewcommand{\d}{\delta}
\newcommand{\D}{\Delta}
\newcommand{\e}{\varepsilon}
\newcommand{\f}{\varphi}
\newcommand{\s}{\sigma}

\renewcommand{\l}{\lambda}

\newcommand{\cS}{{\mathcal S}}

\newcommand{\cV}{{\mathcal V}}

\newcommand{\cR}{{\mathcal R}}
\newcommand{\cK}{{\mathcal K}}

\newcommand{\be}{\begin{equation}}
\newcommand{\ee}{\end{equation}}

\newcommand{\bel}[1]{\begin{equation}\label{#1}}
\newcommand{\beaa}{\begin{eqnarray*}}
\newcommand{\bea}{\begin{eqnarray}}
\newcommand{\beal}[1]{\begin{eqnarray}\label{#1}}
\newcommand{\bean}{\begin{eqnarray}\nonumber}
\newcommand{\beadl}[1]{\begin{deqarr}\label{#1}}
\newcommand{\eeadl}[1]{\arrlabel{#1}\end{deqarr}}
\newcommand{\eeal}[1]{\label{#1}\end{eqnarray}}
\newcommand{\eead}[1]{\end{deqarr}}
\newcommand{\eea}{\end{eqnarray}}
\newcommand{\eeaa}{\end{eqnarray*}}

\renewcommand{\phi}{\varphi}
\renewcommand{\epsilon}{\varepsilon}
\renewcommand{\hat}{\widehat}

\newcommand{\<}{\langle}
\renewcommand{\>}{\rangle}

\newcommand{\dm}{{\partial M}}

\newcommand{\w}{\widetilde}

\theoremstyle{plain}
\newtheorem{theorem}{Theorem}[section]

\newtheorem{Example}[theorem]{Example}
\newtheorem{remark}[theorem]{Remark}

\newtheorem{lemma}[theorem]{Lemma}

\newtheorem{proposition}[theorem]{Proposition}
\newtheorem{corollary}[theorem]{Corollary}

\theoremstyle{definition}

\def\endproof{\qed \medskip}
\def\blacksquare{\hbox to .60em {\vrule width .60em height .60em}}

\numberwithin{equation}{section}

\begin{document}

\title[]{Extension of Symmetries on Einstein Manifolds \\ with Boundary}

\author[]{Michael T. Anderson}

\thanks{Partially supported by NSF Grants DMS 0604735 and DMS 0905159.}

\abstract{We investigate the validity of the isometry extension property for 
(Riemannian) Einstein metrics on compact manifolds $M$ with boundary $\partial M$. 
Given a metric $\gamma$ on $\partial M$, this is the issue of whether any Killing 
field $X$ of $(\partial M, \gamma)$ extends to a Killing field of any Einstein metric 
$(M, g)$ bounding $(\partial M, \gamma)$. Under a mild condition on the fundamental 
group, this is proved to be the case at least when $X$ preserves the mean curvature 
of $\partial M$ in $(M, g)$. }
\endabstract

\maketitle

\setcounter{section}{0}

\section{Introduction.}
\setcounter{equation}{0}

  Let $M^{n+1}$ be a compact $(n+1)$-dimensional manifold-with-boundary, and suppose $g$ 
is a (Riemannian) Einstein metric on $M$, so that
\begin{equation}\label{e1.1}
Ric_{g} = \lambda g,
\end{equation}
for some constant $\lambda\in{\mathbb R}$. The metric $g$ induces a Riemannian boundary 
metric $\gamma$ on $\partial M$. In this paper we consider the issue of whether 
isometries of the boundary structure $(\partial M, \gamma)$ necessarily extend 
to isometries of any filling Einstein manifold $(M, g)$. 

  In general, without any assumptions, this isometry extension property will not hold. 
It is false for instance if $\partial M$ is not connected. For example, let $M = S^{3} \setminus 
(B_{1}\cup B_{2})$, where $B_{i}$ are a pair of disjoint round 3-balls in $S^{3}$ endowed 
with a round metric; then a generic pair of Killing fields $X_{i}$ on $S_{i}^{2} = \partial B_{i}$ 
does not extend to a Killing field on $M$. Also, setting $M = T^{3}\setminus B$ where $B$ is 
a round 3-ball in a flat 3-torus $T^{3}$, one sees again that Killing fields on $\partial M$ 
do not extend to Killing fields on $T^{3}$. This is due to the fact that $\pi_{1}(\partial M)$ 
does not surject onto $\pi_{1}(M)$. Both situations above can be remedied by making the 
topological assumption
\begin{equation}\label{e1.2}
\pi_{1}(M, \partial M) = 0,
\end{equation}
so we will usually assume \eqref{e1.2}. 

  However, this condition is still not sufficient. Consider for example the flat product 
metric on $S^{1}\times {\mathbb R}^{2}$. Let $\sigma$ be any simple closed curve in 
${\mathbb R}^{2}$ and let $T_{\sigma} = S^{1}\times \sigma \subset S^{1}\times 
{\mathbb R}^{2}$. Then $T_{\sigma}$ bounds a compact domain $M \subset S^{1}\times 
{\mathbb R}^{2}$, diffeomorphic to a solid torus. Any such $T_{\sigma}$ is flat 
with respect to the induced metric, and so has a pair of orthogonal Killing fields. 
One of these, that tangent to the $S^{1}$ factor, clearly extends to a Killing field 
of $M$ (in fact $S^{1}\times {\mathbb R}^{2}$). However, whenever $\sigma$ is not a 
round circle in ${\mathbb R}^{2}$ (so that $\sigma$ has non-constant geodesic curvature) 
the orthogonal Killing field on $(T_{\sigma}, \gamma)$ tangent to $\sigma$ does not extend 
as a Killing field to $M$. 

  Very similar examples are easily constructed via the Hopf fibration in the sphere 
${\mathbb S}^{3}$, with $M$ again a solid torus in $S^{3}$, as first pointed out 
to the author by H. Rosenberg \cite{R}, cf.~\cite{Bo}, \cite{Ki}, \cite{Pi} and 
references therein for detailed discussion. Similar examples, even with convex 
boundary, also occur in hyperbolic space-forms, cf.~Remark 3.6 below, and in 
higher dimensions by taking products. 

  The main result of this paper characterizes one situation where the isometry extension 
property does hold. Let $H$ denote the mean curvature of $\partial M$ in $(M, g)$. 

\begin{theorem}\label{t1.1}
Let $g$ be a $C^{m,\alpha}$ Einstein metric on $M$, $m \geq 5$, with induced boundary 
metric $\gamma$ on $\partial M$, and suppose \eqref{e1.2} holds. Then any Killing 
field $X$ on $(\partial M, \gamma)$ for which $X(H) = 0$, extends uniquely to a 
Killing field on $(M, g)$. 
\end{theorem}

   It follows for instance that for $H = const$, the identity component 
$Isom_{0}(\partial M, \gamma)$ of the isometry group of $(\partial M, \gamma)$ 
embeds in the isometry group of any Einstein filling metric $(M, g)$:
$$Isom_{0}(\partial M, \gamma) \hookrightarrow Isom_{0}(M, g),$$
or equivalently, such isometries of the boundary extend to isometries of 
any Einstein filling metric. A simple consequence of Theorem 1.1 is for example 
the following rigidity result.

\begin{corollary}\label{c1.2}
Let $g$ be a $C^{5,\alpha}$ Einstein metric on $M^{n+1}$ which induces the round 
metric $\gamma_{+1}$ on the boundary $\partial M = S^{n}$, $n \geq 2$. If $\pi_{1}(M) 
= 0$ and $H = const$, then $(M, g)$ is isometric to a standard round ball in a 
simply connected space form. 
\end{corollary}

  There are natural analogs of these results valid for exterior domains. Thus, let 
$M^{n+1}$ be an open or non-compact manifold with compact ``inner'' boundary and 
with a finite number of non-compact ends. Metrically, consider complete metrics $g$ 
on $M$ which are asymptotically (locally) flat on each end. In this context, 
Theorem 1.1 also holds for Einstein metrics, cf.~Proposition 5.3. A similar result 
also holds for complete, asymptotically hyperbolic Einstein metrics, with boundary 
at infinity, without any assumption on the mean curvature, cf.~Theorem 5.4. 

\medskip

   We point out that Theorem 1.1 (and Corollary 1.2) remain valid without the hypothesis 
\eqref{e1.2} provided $(M, g)$ is embedded as a domain in a complete, simply connected 
Einstein manifold $(\hat M, \hat g)$. It should also be noted that the isometry extension 
property is false for isometries not contained in $Isom_{0}(\partial M, \gamma)$. As a 
simple example, consider a flat metric on a solid torus $M = D^{2}\times S^{1}$ of the form
$$g_{0} = dr^{2} + r^{2}d\theta_{1}^{2} + d\theta_{2}^{2},$$
for $r \in [0,1]$. Then interchanging the two circles parametrized by $\theta_{1}$ and 
$\theta_{2}$ is an isometry of the boundary, which does not extend to an isometry of 
the solid torus. Of course $\partial M$ is both convex and has constant mean curvature 
in $(M, g_{0})$. 

\medskip

  The proofs of the results above follow from a study of the global properties 
of the space of Einstein metrics $g$ on $M$. As shown in \cite{An1}, the moduli 
space ${\mathcal E}$ of such metrics is a smooth Banach manifold, for which the 
(Dirichlet) map to the boundary metrics
\begin{equation}\label{e1.3}
\Pi_{D}: {\mathcal E} \rightarrow Met(\partial M), \ \ \Pi_{D}(g) = g_{T(\partial M)},
\end{equation}
is $C^{\infty}$ smooth, cf.~Theorem 2.1. The main results are then quite 
simple to prove when the metric $(M, g)$ is non-degenerate, in the strong 
sense that the derivative $D\Pi_{D}$ of $\Pi_{D}$ at $g$ has trivial kernel, 
cf.~Remark 3.3. They also hold, with somewhat more involved proofs, when 
$D\Pi_{D}$ has no cokernel, or more precisely when $Im D\Pi_{D}$ is dense in 
$TMet(\partial M)$, cf.~Proposition 4.4. As discussed in Section 3, note however 
that the map $\Pi_{D}$ is never Fredholm, and the image of the linearization 
is always of infinite codimension. In both of the situations above, the results 
hold without any condition on the mean curvature, i.e.~without assuming $X(H) = 0$. 

  In general, the strategy is to prove the implication
\begin{equation}\label{e1.4}
X(H) = {\mathcal L}_{X}H = 0 \Rightarrow {\mathcal L}_{X}A = 0,
\end{equation}
where $A$ is the $2^{\rm nd}$ fundamental form of $\partial M$ in $M$ and 
${\mathcal L}_{X}$ is the Lie derivative with respect to $X$. Given this, 
Theorem 1.1 then follows from a unique continuation theorem for Einstein metrics 
proved in \cite{AH}. A key point is to relate \eqref{e1.4} with the linearization 
of the divergence constraint for the Einstein equations at $\partial M$, which reads:
\begin{equation}\label{e1.5}
\delta'_{h}(A - H\gamma) + \delta(A - H\gamma)'_{h} = -(Ric(N, \cdot))'_{h} 
\ \ {\rm at} \ \  \partial M,
\end{equation}
where $N$ is the unit normal and $\delta$ is the divergence operator. We show 
in Section 4 (cf.~Lemma 4.2) that \eqref{e1.4} holds provided the linearized divergence 
constraint for the Einstein equations is ``surjective'' at $\partial M$. This 
means that any symmetric form $h^{T}$ on $\partial M$ has an extension to a 
neighborhood of $\partial M$ in $M$ such 
that the derivative 
\begin{equation}\label{e1.6}
(Ric'_{h})(N, \cdot) = 0 \ \ {\rm at} \ \  \partial M,
\end{equation}
This is of course closely related to the surjectivity of $D\Pi_{D}$. Now while 
\eqref{e1.6} does not hold in general (i.e.~for all $h^{T}$ on $\partial M$) 
we prove that any $h$ as above always has an extension such that
\begin{equation}\label{e1.7}
\int_{\partial M}Ric'_{h}(N, \cdot)dV_{\gamma} = 0,
\end{equation}
and this suffices to establish \eqref{e1.4}. The proof of \eqref{e1.7} requires 
a careful study of the linearized Einstein operator, and related operators, with 
certain self-adjoint, elliptic boundary conditions distinct (of course) from 
Dirichlet boundary data; see in particular the operators and boundary conditions 
in \eqref{e5.1} and \eqref{e5.15}. 

\medskip

  A brief survey of the contents of the paper is as follows. In Section 2, we introduce 
the basic setting and structural results on the space of Einstein metrics, needed 
for the work to follow. Section 3 studies elliptic boundary value problems for the 
Einstein equations and the lack of the Fredholm property for the boundary map 
$\Pi_{D}$ in \eqref{e1.3}. Section 4 relates the isometry extension property with 
the linearized divergence constraint equations induced by the Einstein equations 
on $\partial M$. In Section 5, we prove Theorem 1.1 and Corollary 1.2, and the further 
related results mentioned above. 

\medskip

  I would like to thank Mohammad Ghomi and Harold Rosenberg for providing 
very useful background information and references on topics related to this 
paper. Thanks also to one of the referees for suggestions leading to improvements 
in the exposition. Initial work on this paper was carried out at the Institut 
Mittag-Leffler, Djursholm, Sweden, to whom I am grateful for hospitality and 
financial support.

\section{The Space of Einstein Metrics}
\setcounter{equation}{0}

  As above, let $M$ denote a connected, compact, oriented $(n+1)$-dimensional 
manifold with compact, non-empty boundary $\partial M$. Consider the Banach space 
\begin{equation} \label{e2.1}
Met(M) = Met^{m,\alpha}(M) 
\end{equation}
of Riemannian metrics on $M$ which are $C^{m,\alpha}$ smooth up to $\partial M$. 
Here $m$ is any fixed integer with $m \geq 2$, including $m = \infty$ (giving a 
Fr\'echet space) and $\alpha \in (0,1)$. Let 
\begin{equation} \label{e2.2}
{\mathbb E} = {\mathbb E}^{m,\alpha}(M) \subset Met^{m,\alpha}(M)
\end{equation}
be the subset of Einstein metrics on $M$, $C^{m,\alpha}$ smooth up to $\partial M$, 
with 
\begin{equation} \label{e2.3}
Ric_{g} = \lambda g, 
\end{equation}
for $\lambda$ arbitrary, but fixed (so that ${\mathbb E} = {\mathbb E}(\lambda)$); 
$Ric_{g}$ is the Ricci curvature of $g$. The smoothness index $(m, \alpha)$ will 
occasionally be suppressed from the notation when its exact value is unimportant. 

   The space ${\mathbb E}^{m,\alpha}(M) \subset Met^{m,\alpha}(M)$ is invariant under 
the action of the group ${\mathcal D}_{1} = {\mathcal D}_{1}^{m+1,\alpha}$ of orientation 
preserving $C^{m+1,\alpha}$ diffeomorphisms of $M$ equal to the identity on $\partial M$. 
This action is free (since any such isometry equal to the identity on $\partial M$ is 
necessarily the identity) and well-known to be proper. The moduli space ${\mathcal E} 
= {\mathcal E}^{m,\alpha}(M)$ of Einstein metrics on $M$ is defined to be the quotient
\begin{equation} \label{e2.4}
{\mathcal E}  = {\mathbb E}/{\mathcal D}_{1}. 
\end{equation}

   One has a natural Dirichlet boundary map
\begin{equation}\label{e2.5}
\Pi_{D}: {\mathbb E} \rightarrow Met(\partial M); \ \Pi_{D}(g) = \gamma = g|_{T(\partial M)}.
\end{equation}
which clearly descends to a map
\begin{equation} \label{e2.6}
\Pi_{D} : {\mathcal E}  \rightarrow  Met(\partial M);  \ \ \Pi_{D}([g]) = \gamma . 
\end{equation}

  We note the following result, proved in \cite{An1}. 

\begin{theorem}\label{t2.1}
Suppose $\pi_{1}(M, \partial M) = 0$ and $m \geq 5$. Then the space ${\mathcal E}$ is 
a $C^{\infty}$ smooth Banach manifold, {\rm (}Fr\'echet manifold when $m = \infty${\rm )}, 
and the boundary map $\Pi_{D}$ is $C^{\infty}$ smooth. 
\end{theorem}

  Theorem 2.1 is proved by a suitable application of the implicit function 
theorem. Strictly speaking, this result is not needed for the proof of the main results 
in the Introduction; however it places the arguments to follow in a natural context. 

  Consider the Einstein operator
\begin{equation} \label{e2.7}
E: Met(M) \rightarrow  S_{2}(M), 
\end{equation}
$$E(g) = Ric_{g} - \lambda g,$$
where $S_{2}(M)$ is the space of symmetric bilinear forms on $M$. The linearization 
of $E$ is given by
\begin{equation}\label{e2.8}
L_{E}(k) = 2\frac{d}{dt}(Ric_{g+tk} - \lambda (g+tk))|_{t = 0} = \nabla^{*}\nabla k - 
2R(k) - 2\delta^{*}\beta(k);
\end{equation}
here $\delta^{*}X = \frac{1}{2}{\mathcal L}_{X}g$, $\beta (k) = \delta(k) + 
\frac{1}{2}dtr k$ is the Bianchi operator with respect to $g$, $\nabla^{*}\nabla$ is 
the rough Laplacian ($\nabla^{*}\nabla = -\nabla_{e_{i}}\nabla_{e_{i}}$)  and $R(h)$ 
is the action of the curvature tensor on symmetric bilinear forms $k$, cf.~\cite{Be} 
for instance. 

  The tangent space $T_{g}{\mathbb E}$ is given by $Ker L_{E}$. The derivative of the 
Dirichlet boundary map $\Pi_{D}$ in \eqref{e2.5} acts on forms $k$ satisfying $L_{E}(k) 
= 0$ and is given by 
\begin{equation}\label{e2.9}
(D\Pi_{D})_{g}(k) = k^{T}|_{\partial M},
\end{equation}
where $k^{T}$ is the tangential projection or restriction of $k$ to $T(\partial M)$. 
Thus $k^{T}$ is the variation of the boundary metric $\gamma = \Pi_{D}(g)$. It will also 
be important to consider the variation of the $2^{\rm nd}$ fundamental form $A$ of 
$\partial M$ in $M$. Thus, analogous to \eqref{e2.6}, one has a natural Neumann 
boundary map
\begin{equation}\label{e2.10}
\Pi_{N}: {\mathcal E} \rightarrow S^{2}(\partial M), \ \ \Pi_{N}([g]) = A.
\end{equation}
This is well-defined, since $A$ is invariant under the action of ${\mathcal D}_{1}$. 
Note also that $\Pi_{N}$ maps ${\mathcal E}^{m,\alpha}$ to $S_{2}^{m-1,\alpha}
(\partial M)$. To compute the derivative of $\Pi_{N}$, let $g_{s} = g + sk$ be 
a variation of $g$. Since $A = \frac{1}{2}{\mathcal L}_{N}g$, one has $2A'_{k} 
\equiv 2\frac{d}{ds}A_{g_{s}}|_{s=0} =  ({\mathcal L}_{N_{s}}g_{s})'|_{s=0} = 
{\mathcal L}_{N}k + {\mathcal L}_{N'}g$. A simple computation gives $N' = 
- k(N)^{T} - \frac{1}{2}k_{00}N$, where $k(N)^{T}$ is the component of 
$k(N)$ tangent to $\partial M$ and $k_{00} = k(N,N)$. Thus
\begin{equation}\label{e2.11}
A_{k}' = (D\Pi_{N})(k) = \frac{1}{2}({\mathcal L}_{N}k + \delta^{*}V),
\end{equation}
where $V = 2N' = -2k(N)^{T} - k_{00}N$.

  The kernel of $D\Pi_{D}$ in \eqref{e2.5} consists of forms $k$ satisfying 
$L_{E}(k) = 0$ and $k^{T} = 0$ on $\partial M$, while the kernel of $D\Pi_{N}$ in 
\eqref{e2.10} consists of such forms satisfying $(A_{k}')^{T} = 0$ at $\partial M$. 
Thus, if both conditions hold, 
\begin{equation}\label{e2.12}
k^{T} = 0, \ (A_{k}')^{T} = 0 \ \ {\rm at} \ \ {\partial M},
\end{equation}
then $(M, g)$ is both Dirichlet and Neumann degenerate, i.e.~a singular point of each 
boundary map. We note that each of the conditions in \eqref{e2.12} is gauge-invariant, 
i.e.~invariant under the addition of terms of the form $\delta^{*}Z$ with $Z = 0$ 
on $\partial M$. Of course any form $k$ satisfying $k = \nabla_{N}k = 0$ at 
$\partial M$ satisfies \eqref{e2.12}. Changing such $k$ by arbitrary such gauge 
transformations shows that \eqref{e2.12} is equivalent to the statement that $k$ 
is pure gauge, to first order at $\partial M$, i.e.
\begin{equation}\label{e2.13}
k = \delta^{*}Z + O(t^{2}),
\end{equation}
near $\partial M$, with $Z = 0$ on $\partial M$, where $t(x) = dist_{g}(x, \partial M)$. 

  The natural or geometric Cauchy data for the Einstein equations \eqref{e2.3} on $M$ at 
$\partial M$ consist of the pair $(\gamma, A)$. If $k$ is an infinitesimal Einstein 
deformation of $(M, g)$, so that $L_{E}(k) = 0$, then the induced variation of the 
Cauchy data on $\partial M$ is given by $k^{T}$ and $(A_{k}')^{T}$. It is natural 
to expect that an Einstein metric $g$ is uniquely determined in a neighborhood of 
$\partial M$, up to isometry, by the Cauchy data $(\gamma, A)$, i.e.~one should have 
a suitable unique continuation property for Einstein metrics. Similarly, one would 
expect this holds for the linearized Einstein equations. The next result, proved in 
\cite{AH} confirms this expectation; this result is also proved in \cite{Bi}, but 
only under the much stronger assumption of $C^{\infty}$ smoothness up to the 
boundary.

\begin{theorem}\label{et2.2}
Let $g \in {\mathbb E}^{m,\alpha}$, $m \geq 5$, and suppose $k$ is an infinitesimal 
Einstein deformation which is both Dirichlet and Neumann degenerate, so that 
$L_{E}(k) = 0$ and \eqref{e2.12} holds. Then $k$ is pure gauge near $\partial M$, i.e. 
\begin{equation}\label{e2.14}
k = \delta^{*}Z \ \ {\rm near} \ \ \partial M, 
\end{equation}
with $Z = 0$ on $\partial M$. 
\end{theorem}

  As is well-known, the operator $E$ is not elliptic, due to its covariance under 
diffeomorphisms: one has $L_{E}(\delta^{*}Y) = 0$, for any vector field $Y$ on $M$, 
at an Einstein metric. We will require ellipticity at several points and so need 
a choice of gauge to break the diffeomorphism invariance of the Einstein equations. 
In view of \eqref{e2.8}, the simplest and most natural choice is the Bianchi gauge 
given by $\beta (k) = 0$ at the linearized level. (Later we will use instead a 
slightly different gauge, the divergence-free gauge, cf.~Remark 4.6). Thus, 
let $\widetilde g$ be a fixed (background) metric in ${\mathbb E}$. The 
associated Bianchi-gauged Einstein operator is given by the $C^{\infty}$ smooth 
map
\begin{equation} \label{e2.15}
\Phi_{\widetilde g}: Met^{m,\alpha}(M) \rightarrow  S_{2}^{m-2,\alpha}(M), 
\end{equation}
$$\Phi(g) = \Phi_{\widetilde g}(g) = Ric_{g} - \lambda g + \delta_{g}^{*}
\beta_{\widetilde g}(g),$$
where $\beta_{\widetilde g}(g)$ is the Bianchi operator with respect to 
$\widetilde g$, while $\delta^{*}$ is taken with respect to $g$. Although 
$\Phi_{\widetilde g}$ is defined for all $g\in Met(M)$, we will only consider 
it acting on $g$ near $\widetilde g$. 

  The linearization 
of $\Phi$ at $\widetilde g = g$ is given by 
\begin{equation} \label{e2.16}
L: T_{\widetilde g}Met(M) \rightarrow S_{2}(M),
\end{equation}
$$L(h) = 2(D\Phi)_{\widetilde g}(h) = \nabla^{*}\nabla h - 2R(h).$$
The operator $L$ is formally self-adjoint and is clearly elliptic. Comparing 
\eqref{e2.7} and \eqref{e2.15}, the relation between $L$ and the linearization 
$L_{E} = 2E'$ of the Einstein operator $E$ in \eqref{e2.8} is given by 
\begin{equation}\label{e2.17}
L_{E} = L - 2\delta^{*}\beta .
\end{equation}
In Section 3, we will consider elliptic boundary value problems for the operator $\Phi$. 

 Clearly $g \in {\mathbb E}$ if $\Phi_{\widetilde g}(g) = 0$ and 
$\beta_{\widetilde g}(g) = 0$, so that $g$ is in the Bianchi gauge with 
respect to $\widetilde g$. Given $\widetilde g$, let $Met_{C}(M) 
= Met_{C}^{m,\alpha}(M)$ be the space of $C^{m,\alpha}$ smooth Riemannian metrics 
on $M$ in Bianchi gauge with respect to $\widetilde g$ at $\partial M$:
\begin{equation}\label{e2.18}
Met_{C}(M) = \{g \in Met(M):\beta_{\widetilde g}(g) = 0 \ {\rm at} \ \partial M\}.
\end{equation}
Let 
\begin{equation}\label{e2.19}
Z_{C} = \{g \in Met_{C}(M): \Phi(g) = 0 \}
\end{equation}
be the $0$-set of $\Phi$ and let ${\mathbb E}_{C} \subset Z_{C}$ be the 
subset of Einstein metrics $g$ in $Z_{C}$. 
 
  To justify the use of $\Phi$, one needs to show that the opposite inclusion holds, 
so that ${\mathbb E}_{C} = Z_{C}$. This has already been done in \cite{An1} and we 
summarize the results here. 
\begin{lemma}\label{l2.3}
{\bf (i).} For $g$ in $Met^{m,\alpha}(M)$, one has
\begin{equation}\label{e2.20}
T_{g}Met^{m-2,\alpha}(M) \simeq S_{2}^{m-2,\alpha}(M) = Ker \delta \oplus 
Im \delta^{*},
\end{equation}
where $\delta^{*}$ acts on $\chi_{1}^{m-1,\alpha}$, the space of $C^{m-1,\alpha}$ 
vector fields on $M$ which vanish on $\partial M$. 

{\bf (ii).} For $\widetilde g \in {\mathbb E}^{m,\alpha}$ and $g$ in 
$Met^{m,\alpha}(M)$ close to $\widetilde g$, one has
\begin{equation}\label{e2.21}
T_{g}Met^{m-2,\alpha}(M) \simeq S_{2}^{m-2,\alpha}(M) = Ker \beta \oplus 
Im \delta^{*},
\end{equation}
where $\beta$ is the Bianchi operator with respect to $\widetilde g$. If 
$g \in {\mathbb E}^{m,\alpha}$, then \eqref{e2.21} holds with $m$ in place 
of $m-2$. 

  {\bf (iii).} 
Any metric $g \in Z_{C}$ near $\widetilde g \in {\mathbb E}^{m,\alpha}$ is 
Einstein, and in Bianchi gauge with respect to $\widetilde g$, i.e.
\begin{equation}\label{e2.22}
\beta_{\widetilde g}(g) = 0.
\end{equation}
Similarly, if $k \in Met_{C}(M)$ is an infinitesimal deformation of $\widetilde g$ 
in $Z_{C}$, i.e.~$L(k) = 0$, then $k$ is an infinitesimal Einstein deformation 
and $\beta (k) = 0$. 
\end{lemma}

   Lemma 2.3 implies that ${\mathbb E}_{C} = Z_{C}$ near $\widetilde g$, and at 
least infinitesimally ${\mathbb E}_{C}$ is a local slice for the action of the 
diffeomorphism group ${\mathcal D}_{1}$ on ${\mathbb E}$. In fact, it is shown 
in \cite{An1} that ${\mathbb E}_{C}$ is a local slice for the action of 
${\mathcal D}_{1}$. 

  The next two results may be viewed as a preliminary version of Theorem 1.1.

\begin{corollary}\label{c2.4}
Let $g \in {\mathbb E}^{m,\alpha}$, $m \geq 5$, and suppose $\kappa$ is an 
infinitesimal Einstein deformation of $(M, g)$. If $\pi_{1}(M, \partial M) = 0$ 
and \eqref{e2.12} holds, then $\kappa$ is pure gauge on $M$, i.e.~there exists 
a vector field $Z$ on $M$ with $Z = 0$ on $\partial M$ such that
\begin{equation}\label{e2.23}
\kappa = \delta^{*}Z \ \ {\rm on} \ \ M.
\end{equation}
If $\kappa$ is in Bianchi gauge, so that $L(\kappa) = 0$, then 
$$\kappa = 0 \ \ {\rm on} \ \  M.$$ 
\end{corollary}

{\bf Proof:}  The hypotheses and Theorem 2.2 imply that the form $\kappa$ on 
$M$ is pure gauge near $\partial M$, so that \eqref{e2.23} holds on a 
neighborhood $\Omega$ of $\partial M$. 

  It then follows from a basically standard analytic continuation argument in the 
interior of $M$ that the vector field $Z$ may be extended so that \eqref{e2.23} 
holds on all of $M$, cf.~[11, Ch.~VI.6.3] for instance. A detailed proof of this is 
also given in [3, Lemma 2.6]. This analytic continuation argument requires the 
topological hypothesis \eqref{e1.2} to obtain a well-defined (single-valued) 
vector field $Z$ on $M$. Moreover, since $\partial M$ is connected, the condition 
$Z = 0$ on $\partial M$ remains valid in the analytic continuation. 

  For the second statement, if in addition $\beta (\kappa) = 0$, then 
$\beta \delta^{*}Z = 0$ on $M$ with $Z = 0$ on $\partial M$. It then follows 
from Lemma 2.3 that $Z = 0$ on $M$ and hence $\kappa = 0$ on $M$, as claimed. 

{\endproof}

\begin{proposition}\label{p2.5}
Let $g \in {\mathbb E}^{m,\alpha}$, $m \geq 5$, and suppose $X$ is a Killing 
field on $(\partial M, \gamma)$ such that
\begin{equation}\label{e2.24}
({\mathcal L}_{X}A)^{T} = 0 \ \ {\rm at} \ \ \partial M,
\end{equation}
If $\pi_{1}(M, \partial M) = 0$, then $X$ extends to a Killing field on $(M, g)$. 
\end{proposition}

{\bf Proof:} Since $\gamma \in Met^{m,\alpha}(\partial M)$, the Killing field $X$ is 
$C^{m+1,\alpha}$ smooth on $\partial M$. By Lemma 2.3, the operator $\beta \delta^{*}$ 
has no kernel on $\chi_{1}$. Since this operator has (Fredholm) index 0, it 
follows that $X$ may be uniquely extended to a vector field $X$ on $M$ so that
\begin{equation}\label{e2.25}
\beta \delta^{*}X = 0 \ \ {\rm on} \ \ M.
\end{equation}
Since $g \in {\mathbb E}^{m,\alpha}$, the solution $X$ is then $C^{m+1,\alpha}$ up to 
$\partial M$. Hence the form $\kappa = \delta^{*}X$ is $C^{m-1,\alpha}$ up to $\partial M$ 
and is an infinitesimal Einstein deformation in Bianchi gauge, i.e.~$L(\kappa) = 
L_{E}(\kappa) = 0$ with $\beta (\kappa) = 0$. Note that by construction, $\kappa^{T} 
= 0$ at $\partial M$. 

  Next, note that 
\begin{equation}\label{e2.26}
{\mathcal L}_{X}A = 2A'_{\kappa}.
\end{equation}
Namely, since $\kappa = \frac{1}{2}{\mathcal L}_{X}g$, as in \eqref{e2.11} one 
has $A'_{\kappa} = \frac{1}{4}{\mathcal L}_{N}{\mathcal L}_{X}g + \frac{1}{2}
{\mathcal L}_{N'}g = \frac{1}{2}{\mathcal L}_{X}A + \frac{1}{4}{\mathcal L}_{[N,X]}g 
+ \frac{1}{2}{\mathcal L}_{N'}g$. It is easy to verify that $[N, X] = -2N'$, which 
gives \eqref{e2.26}. 

  The form $\kappa$ is thus an infinitesimal Einstein deformation in Bianchi gauge 
and satisfies \eqref{e2.12}. Hence by Corollary 2.4, $\kappa = \delta^{*}X = 0$ 
on $M$, which implies that $X$ is a Killing field on $(M, g)$. 

{\endproof}

  We note that if, in place of the condition $\pi_{1}(M, \partial M) = 0$, one assumes 
that $(M, g)$ is embedded as a domain in a complete, simply connected Einstein manifold 
$(\hat M, \hat g)$, then essentially the same analytic continuation argument 
(cf.~[9, Ch.~VI.6.4]) again implies that $X$ extends uniquely to a Killing field 
on all of $\hat M$, which proves Proposition 2.5 in this case also.

\section{Elliptic Boundary Value Problems for the Einstein Equations}
\setcounter{equation}{0}

  In this section, we consider elliptic boundary value problems for the Bianchi-gauged 
Einstein operator $\Phi$ in \eqref{e2.15} and the Fredholm properties of the Dirichlet 
boundary map $\Pi_{D}$ in \eqref{e2.6}. 

  Recall that the kernel of the linearized operator $L$ in \eqref{e2.16} forms the 
tangent space $T_{g}Z_{C}$ ($g = \widetilde g$ here) and by Lemma 2.3,
\begin{equation}\label{e3.1}
T_{g}Z_{C} = T_{g}{\mathbb E}_{C},
\end{equation}
so that the kernel also represents the space of (non-trivial) infinitesimal Einstein 
deformations in Bianchi gauge. The natural Dirichlet-type boundary conditions for 
$\Phi$ are
\begin{equation}\label{e3.2}
\beta_{\widetilde g} (g) = 0, \ \ g^{T} = \gamma \ \ {\rm at} \ \ \partial M.
\end{equation}
However, contrary to first impressions, the operator $\Phi$ with boundary conditions 
\eqref{e3.2} does not form a well-defined elliptic boundary value problem (for $g$ 
near $\widetilde g$). This is due to the well-known constraint equations, induced 
by the Gauss and Gauss-Codazzi equations on $\partial M$:
\begin{equation} \label{e3.3}
\delta(A - H\gamma) = -Ric_{g}(N, \cdot) = 0,
\end{equation}
\begin{equation} \label{e3.4}
|A|^{2} - H^{2} + s_{\gamma} = s_{g} - 2Ric_{g}(N,N) = (n-1)\lambda.
\end{equation}
Here $H$ is the mean curvature of $\partial M$ in $M$, while $s$ denotes the 
scalar curvature. 

  As will be seen in Section 4, the momentum or vector constraint \eqref{e3.3} is an 
important issue in the study of the isometry extension or rigidity results discussed 
in the Introduction. On the other hand, the Hamiltonian or scalar constraint \eqref{e3.4} 
is important in understanding the Fredholm properties of the boundary map $\Pi_{D}$ in 
\eqref{e2.6}. Thus for $g \in {\mathbb E}^{m,\alpha}$, one has $A \in S_{2}^{m-1,\alpha}
(\partial M)$ so that \eqref{e3.4} implies that $s_{\gamma} \in C^{m-1,\alpha}(\partial M)$. 
However, the space of metrics $\gamma \in Met^{m,\alpha}(\partial M)$ for which $s_{\gamma} 
\in C^{m-1,\alpha}(\partial M)$ is of infinite codimension in $Met^{m,\alpha}(\partial M)$.  
It follows that the linearization of the boundary map $\Pi_{D}$ has infinite dimensional cokernel, 
at least when $m < \infty$, and so $\Pi_{D}$ is never Fredholm. Hence, the boundary 
conditions \eqref{e3.2} for the operator $\Phi$ are not elliptic.

\begin{remark}\label{r3.1}
{\rm It is worthwhile to understand situations where the linearization $D\Pi_{D}$ has 
infinite dimensional kernel and cokernel, even in the $C^{\infty}$ case. Let
\begin{equation}\label{e3.5}
K = K_{g} = Ker D_{g}\Pi_{D}.
\end{equation}
Via the slice representation $Z_{C} = {\mathbb E}_{C} \subset {\mathbb E}$ at 
$\widetilde g = g$, $K$ consists of forms $\kappa$ such that 
\begin{equation}\label{e3.6}
L(\kappa) = 0 \ {\rm and} \ \beta_{g}(\kappa) = 0 \ {\rm on} \  M, 
\ {\rm with} \ \kappa^{T} = 0 \ {\rm on} \ \partial M.
\end{equation}

  Consider then the intersection $K \cap Im \delta^{*}$. Let $Y$ be a vector field 
at $\partial M$ (not necessarily tangent to $\partial M$) and extend $Y$ to 
a vector field on $M$ to be the unique solution to the equation $\beta(\delta^{*}Y) 
= 0$ with the given boundary value, cf.~Lemma 2.3. Then $L(\delta^{*}Y) = 0$ and the 
boundary condition $k^{T} = (\delta^{*}Y)^{T} = 0$ is equivalent to the equation
\begin{equation}\label{e3.7}
(\delta^{*}Y^{T})^{T} + \langle Y, N \rangle A = 0 \ \ {\rm at} \ \ \partial M.
\end{equation}
In particular if $\delta^{*}_{T}$ is the restriction of $\delta^{*}$ to vector 
fields tangent to $\partial M$ at $\partial M$, then $K \cap Im \delta^{*}_{T}$ 
is isomorphic to the space of Killing fields on $(\partial M, \gamma)$. 

  On the other hand, if $\partial M$ is totally geodesic on some open set 
$U \subset \partial M$, i.e.~$A = 0$ on $U$, then the system \eqref{e3.7} has 
solutions of the form $Y = fN$, for {\it any} $f$ with $supp \,f \subset U$, so 
that $K \cap Im \delta^{*}$ is infinite dimensional. Such vector fields $Y$ are 
infinitesimal isometries {\it at} (as opposed to {\it on}) $\partial M$, in 
that they preserve the metric $\gamma$ on $\partial M$ to first order. Of course 
in general such $Y$ do not extend to a Killing field on $(M, g)$; see also Remark 4.3 
for further discussion and examples. This behavior is classically very well-known 
in the context of surfaces embedded in ${\mathbb R}^{3}$, cf.~\cite{Sp}, \cite{Bo}. 

   A similar phenomenon holds for the cokernel. Thus, suppose $(\partial M, \gamma)$ 
is totally geodesic in a domain $U \subset \partial M$. Consider the linearization 
$s_{\gamma}'(h)$, for $h \in Im(D\Pi_{D})$. By differentiating the scalar constraint 
\eqref{e3.4} in the direction $h$, one sees that $s_{\gamma}'(h) = 0$ on $U$, 
for any such $h$. It follows that $Im D\Pi_{D}$ has infinite codimension, even 
in the $C^{\infty}$ case, in such situations. The same argument and conclusion 
holds if $A = 0$ at just one point in $\partial M$. 

  Very little seems to be understood in characterizing the situations where $K$ 
is finite dimensional or $K = 0$. Again, this is the case even in the classical 
setting of closed surfaces embedded in ${\mathbb R}^{3}$. }
\end{remark}

  The discussion above implies there is no natural elliptic boundary value problem 
for the Einstein equations associated with Dirichlet boundary values. To obtain 
an elliptic problem, one needs to add either gauge-dependent terms or terms 
depending on the extrinsic geometry of $\partial M$ in $(M, g)$. To maintain a 
determined boundary value problem, one then has to subtract part of  the 
intrinsic Dirichlet boundary data on $\partial M$.

  There are several ways to carry this out in practice, but we will concentrate 
on the following situation. First, ellipticity of the Bianchi-gauged Einstein operator 
$\Phi = \Phi_{\widetilde g}$ with respect to given boundary conditions - near a given 
solution - depends only on the linearized operator, so we assume $g = \widetilde g$ is 
Einstein and study the linearized operator $L$ from \eqref{e2.16} at $(M, g)$. As usual, 
let $\gamma$ be the induced metric on $\partial M$. 

   Let $B$ be a $C^{m,\alpha}$ symmetric bilinear form on $\partial M$ such that 
\begin{equation}\label{e3.8}
\tau_{B} = B - (tr_{\gamma}B)\gamma < 0,
\end{equation}
is negative definite; all the statements to follow hold equally well if $\tau_{B}$ is 
positive definite. This condition is equivalent to the statement that the sum of 
any $(n-1)$-eigenvalues of $B$ with respect to $\gamma$ is positive. For the choice 
$B = A$, the $2^{\rm nd}$ fundamental form, this is just the statement $\partial M$ 
is $(n-1)$-convex in $(M, g)$, cf.~\eqref{e3.23} below. 

  In place of prescribing the boundary metric $g^{T}$ or its linearization $h^{T}$ on 
$\partial M$, only $h^{T}$ modulo $B$ will be prescribed. Thus, let $\pi_{B}: 
T_{\gamma}Met^{m,\alpha}(\partial M) \rightarrow S_{2}^{m,\alpha}(\partial M)/B$, 
be the natural projection and set $\pi_{B}(h) = [h^{T}]_{B}$. In place of 
the second equation in \eqref{e3.2}, we impose
\begin{equation}\label{e3.9}
[h^{T}]_{B} = h_{1}.
\end{equation}
For example, when $B$ equals the boundary metric $\gamma$, one is prescribing the 
trace-free part of $h^{T}$, i.e.~the tangent space of conformal classes on $\partial M$. 
Another natural choice is $B = A$, the $2^{\rm nd}$ fundamental form of $\partial M$. 
In this case, for regularity purposes, one must work instead with a smooth 
approximation to $A$, since $A \in S_{2}^{m-1,\alpha}(\partial M)$, or with a 
$C^{\infty}$ background $(M, g)$. 

  The simplest gauge-dependent term one can add to \eqref{e3.2} is the equation 
$h(N, N) = h_{00}$, where $N$ is the unit normal with respect to $g$, while the 
simplest extrinsic geometric scalar is the linearization $H'_{h}$ of the mean 
curvature of $\partial M$ in $(M, g)$ in the direction $h$. As shown in \cite{An1}, 
ellipticity holds for either of these boundary conditions. We will use a slightly 
more general result, whose proof is a simple modification of the proof in 
\cite{An1}. 

\begin{proposition}\label{p3.2}
Suppose $B \in S_{2}^{m,\alpha}$ satisfies \eqref{e3.8} and suppose $\sigma$ 
is any positive definite form in $S_{2}^{m,\alpha}(\partial M)$. Then the 
Bianchi-gauged linearized Einstein operator $L$ in \eqref{e2.16} with boundary 
conditions
\begin{equation}\label{e3.10}
\beta(h) = 0, \ \ [h^{T}]_{B} = h_{1}, \ \langle A'_{h}, \sigma \rangle = 
tr_{\sigma}A'_{h} = h_{2} \ \ {\rm at} \ \ \partial M,
\end{equation}
is an elliptic boundary value problem of Fredholm index 0.
\end{proposition} 

{\bf Proof:} The leading order symbol of $L = D\Phi$ is given by 
\begin{equation}\label{e3.11}
\sigma(L) = -|\xi|^{2}I,
\end{equation}
where $I$ is the $N\times N$ identity matrix, with $N = (n+2)(n+1)/2$ the dimension 
of the space of symmetric bilinear forms on ${\mathbb R}^{n+1}$. In the following, 
the subscript 0 represents the direction normal to $\partial M$ in $M$, and Latin 
indices run from $1$ to $n$. The positive roots of \eqref{e3.11} are $i|\xi|$, with 
multiplicity $N$. 

  Writing $\xi = (z, \xi_{i})$, the symbols of the leading order terms in the boundary 
operators in \eqref{e3.10} are given by:
$$-2izh_{0k} - 2i\sum \xi_{j}h_{jk} + i\xi_{k}tr h = 0,$$
$$-2izh_{00} - 2i\sum \xi_{k}h_{0k} + iztr h = 0,$$
$$h^{T} = h_{1} \ \ mod\, B,$$
$$tr_{\sigma}A'_{h} = h_{2},$$
where $h$ is an $(n+1)\times (n+1)$ matrix. Then ellipticity requires that the operator 
defined by the boundary symbols above has trivial kernel when $z$ is set to the root 
$i|\xi|$. Carrying this out then gives the system 
\begin{equation}\label{e3.12}
2|\xi|h_{0k} - 2i\sum \xi_{j}h_{jk} + i\xi_{k}tr h = 0,
\end{equation}
\begin{equation}\label{e3.13}
2|\xi|h_{00} - 2i\sum \xi_{k}h_{0k} - |\xi|tr h = 0,
\end{equation}
\begin{equation}\label{e3.14}
h_{kl} = \phi b_{kk}\delta_{kl},
\end{equation}
\begin{equation}\label{e3.15}
tr_{\sigma}A'_{h} = 0,
\end{equation}
where without loss of generality we assume $B$ is diagonal, with entries $b_{kk}$, and 
$\phi$ is an undetermined function. 

  Multiplying \eqref{e3.12} by $i\xi_{k}$ and summing gives
$$2|\xi|i\sum \xi_{k}h_{0k} =  2i^{2}\xi_{k}^{2}h_{kk} - i^{2}\xi_{k}^{2}tr h.$$
Substituting \eqref{e3.13} on the term on the left above then gives
$$2|\xi|^{2}h_{00} - |\xi|^{2}tr h =  -2\sum \xi_{k}^{2}h_{kk} + |\xi|^{2}tr h,$$
so that
$$|\xi|^{2}h_{00} - |\xi|^{2}tr h =  -\sum \xi_{k}^{2}h_{kk} = 
-\phi \langle B(\xi), \xi \rangle.$$
Using the fact that $\sum h_{kk} = tr h - h_{00}$, this is equivalent to 
$$\phi\langle B(\xi), \xi \rangle = \phi |\xi|^{2}tr B.$$
Since $\tau_{B} = B - (tr B)\gamma$ is assumed to be definite, it follows 
that $\phi = 0$ and hence $h^{T} = 0$. 

  Next, a simple computation from \eqref{e2.11} shows that to leading order, 
$tr_{\sigma}A'_{h} = tr_{\sigma}(\nabla_{N}h - 2\delta^{*}(h(N)^{T}))$, which has 
symbol $iz\sigma^{ij} h_{ij} - 2i\sigma^{ij}\xi_{i}h_{0j}$. Setting this to 0 at 
the root $z = i|\xi|$ and using the fact that $h^{T} = 0$ gives
\begin{equation}\label{e3.16}
\sigma^{ij}\xi_{i}h_{0j} = 0.
\end{equation}
Now \eqref{e3.12} and $h^{T} = 0$ gives $2|\xi|h_{0j} + i\xi_{j}h_{00} = 0$. Multiplying 
the first term by $\sigma^{ij}\xi_{i}$ and summing over $i,j$ gives 0 by \eqref{e3.16}, 
and hence $\sigma^{ij}\xi_{i}\xi_{j}h_{00} = 0$. Since $\sigma > 0$, it follows that 
$h_{00} = 0$ and hence by \eqref{e3.12} again, $h_{0k} = 0$ for all $k$. This gives 
$h = 0$, and hence the boundary data \eqref{e3.10} are elliptic. 

  To prove the operator $L$ with boundary data \eqref{e3.10} is of Fredholm index 0, 
one may continuously deform the boundary data through elliptic boundary values 
to self-adjoint boundary data, which clearly has index 0. This is done in detail 
for the case $\sigma = \gamma$ in \cite{An1} and the proof for general $\sigma > 0$ 
is identical. Thus we refer to \cite{An1} for details as needed. The result then 
follows from the homotopy invariance of the index. 

{\endproof}

  Given $\widetilde g \in {\mathbb E}^{m,\alpha}$, and $B$ as in \eqref{e3.8}, let 
$Met_{B}^{m,\alpha}(\partial M) = Met^{m,\alpha}(\partial M)/B$ be the space of 
equivalence classes of $C^{m,\alpha}$ metrics on $\partial M$ (mod $B$), with 
natural projection or quotient map 
$$\pi_{B}: Met^{m,\alpha}(\partial M) \rightarrow Met_{B}^{m,\alpha}(\partial M).$$ 
It follows from Proposition 3.2 and Lemma 2.3 that the map
\begin{equation}\label{e3.17}
\widetilde \Pi_{B.\sigma}: {\mathbb E}_{C} \rightarrow Met_{B}^{m,\alpha}(\partial M) 
\times C^{m-1,\alpha}(\partial M),
\end{equation}
$$\widetilde \Pi_{B,\sigma}(g) = ([g^{T}]_{B}, tr_{\sigma}A),$$ 
is Fredholm, of index 0, for $g$ near $\widetilde g$. 

  In analogy to \eqref{e3.5}, let
\begin{equation}\label{e3.18}
\widetilde K_{B,\sigma} = Ker D\widetilde \Pi_{B,\sigma},
\end{equation}
where the derivative is taken at $g = \widetilde g$. In contrast to $K$ in 
\eqref{e3.5}, $\widetilde K_{B,\sigma}$ is always finite dimensional. One might call 
an Einstein metric $g \in {\mathbb E}$ {\it non-degenerate} (or $(B,\sigma)$-nondegenerate) 
if 
\begin{equation} \label{e3.19}
\widetilde K_{B,\sigma} = 0, 
\end{equation}
for some $B$, $\sigma$. Thus, $g$ is non-degenerate if and only if $g$ is a regular 
point of the boundary map $\widetilde \Pi_{B,\sigma}$ in which case 
$\widetilde \Pi_{B,\sigma}$ is a local diffeomorphism near $g$. 

\begin{remark}\label{r3.3}
{\rm It is worth pointing out that if $(M, g)$ is strongly non-degenerate, in the 
sense that $K = 0$ in \eqref{e3.5}, then Theorem 1.1 is easy to prove and holds 
without the assumptions on $H$ or on $\pi_{1}(M, \partial M)$. To see this, let 
$\phi_{s}$ be a local curve of $C^{m+1,\alpha}$ diffeomorphisms of $\bar M$ with 
$\phi_{0} = id$ such that $\frac{d}{ds}\phi_{s}|_{s=0} = X$. If $X$ is a Killing 
field on $(\partial M, \gamma)$, then
$$\phi_{s}^{*}\gamma  = \gamma + O(s^{2}).$$ 
The curve $g_{s} = \phi_{s}^{*}g$ is a smooth curve in ${\mathbb E}$, and by 
construction, one has $[h] = [\frac{dg_{s}}{ds}] \in Ker D\Pi_{D}$, for $\Pi_{D}$ as 
in \eqref{e2.6}. One may then alter the diffeomorphisms $\phi_{s}$ by composition with 
diffeomorphisms $\psi_{s} \in {\mathcal D}_{1}^{m+1,\alpha}$ if necessary, so that 
$\kappa = \frac{d\psi_{s}^{*}(g_{s})}{ds} \in K_{g}$, where $K = K_{g}$ is the kernel 
in \eqref{e3.5} and $[h] = [\kappa]$. Thus
$$\kappa = \delta^{*}X',$$
where $X' = \frac{d(\phi_{s}\circ \psi_{s})}{ds}$ is $C^{m+1,\alpha}$ smooth up to 
$\bar M$. Note that $X' = X$ at $\partial M$. If $K_{g} = 0$, then this gives
$$\delta^{*}X' = 0 \ {\rm on} \ M,$$
so that $X'$ is a Killing field on $(M, g)$. Thus, any Killing field on 
$(\partial M, \gamma)$ extends to a Killing field on $(M, g)$, as claimed. 
The same result and proof hold in general, for any infinitesimal Einstein 
deformation preserving the boundary metric $(\partial M, \gamma)$. 

  It follows that if this general isometry extension property fails, then the 
Dirichlet boundary map $\Pi_{D}$ in \eqref{e3.5} is necessarily degenerate. }
\end{remark}

\begin{remark}\label{r3.4}
{\rm Although currently the cokernel of $D\Pi_{D}$ remains hard to understand, 
cf.~Remark 3.1, it is not difficult to describe the cokernel of 
$D\widetilde \Pi_{B,\sigma}$. For simplicity, set $(B, \sigma) = (\gamma, \gamma)$ 
and let $\widetilde \Pi_{\gamma, \gamma} = \widetilde \Pi_{H}$. Then define
\begin{equation}\label{e3.20}
\widetilde C = \{(({\mathcal L}_{N}\kappa)_{\gamma}^{T}, N(H'_{\kappa})) : \kappa \in 
\widetilde K_{\gamma,\gamma} \},
\end{equation}
so that $\widetilde C$ represents Neumann-type data associated with the Dirichlet 
data in \eqref{e3.9}. 

  Note that $\widetilde C \subset S_{\gamma}^{m,\alpha}(\partial M)\times C^{m-1,\alpha}
(\partial M)$, where $S_{\gamma}^{m,\alpha}(\partial M) = T_{g}Met_{\gamma}^{m,\alpha}
(\partial M) \simeq S^{m,\alpha}(\partial M)/\langle \gamma \rangle$. Namely, for 
$\kappa \in \widetilde K_{\gamma,\gamma}$, one has $L(\kappa) = 0$ on $M$ together 
with the elliptic boundary conditions $\beta (\kappa) = 0$, $\kappa_{\gamma}^{T} = 0$, 
and $H'_{\kappa} = 0$ on $\partial M$. Since $g$ is $C^{m,\alpha}$ up to $\partial M$, 
elliptic boundary regularity applied to this system gives $\kappa \in C^{m+1,\alpha}$ 
(cf.~\cite{GT, M}) so that ${\mathcal L}_{N}\kappa \in S_{2}^{m,\alpha}(\partial M)$ 
and $N(H'_{\kappa}) \in C^{m-1,\alpha}(\partial M)$. 

  It is then not difficult to prove (although we will not give the proof here) 
that the space $\widetilde C$ is a slice for $Coker D\widetilde \Pi_{H}$ in 
$S_{\gamma}^{m,\alpha}(\partial M)\times C^{m,\alpha}(\partial M)$, so that
\begin{equation}\label{e3.21}
S_{\gamma}^{m,\alpha}(\partial M)\times C^{m-1,\alpha}(\partial M) = 
Im D\widetilde \Pi_{H} \oplus \widetilde C .
\end{equation}

  By restricting to the first factor, it follows immediately from \eqref{e3.21} that
\begin{equation}\label{e3.22}
S_{\gamma}^{m,\alpha}(\partial M) = 
Im D\Pi_{0} \oplus \widetilde S,
\end{equation}
where $\widetilde S = \{({\mathcal L}_{N}\kappa)_{\gamma}^{T}: \kappa \in 
\widetilde K_{H}\}$ and $\Pi_{0}$ is defined by $\Pi_{0} = \pi_{\gamma}\circ \Pi_{D}$. }
\end{remark}

   One may use the diffeomorphism group to pass from the space ${\mathbb E}_{C}$ of 
Bianchi-gauged Einstein metrics to the full space ${\mathbb E}$, thus passing from 
$\widetilde \Pi_{H}$ to the more natural Dirichlet boundary map $\Pi_{D}$. In more detail, 
the image ${\mathcal V} = D\Pi_{D}({\mathbb E}_{C}) \subset TMet^{m,\alpha}(\partial M)$ 
projects onto a space  of finite codimension in $S_{\gamma}^{m,\alpha}(\partial M)$ by 
\eqref{e3.22}. The full image $D\Pi_{D}({\mathbb E})$ then consists of the span 
$\langle {\mathcal V}, Im \delta^{*} \rangle$, where $\delta^{*}$ acts on all 
vector fields at $\partial M$, not necessarily tangent to $\partial M$. It is 
an interesting question to understand when the closure of this space is of 
finite codimension in $TMet^{m,\alpha}(\partial M)$. This corresponds roughly 
to $\Pi_{D}$ being Fredholm.

\medskip

   One situation where this occurs is the following. Define $\partial M \subset M$ to 
be $p$-convex if the sum of any $p$ eigenvalues of the second fundamental form $A$ of 
$\partial M$ in $(M, g)$ is positive, cf.~also \cite{Sh} for example. Thus, $\partial M$ 
is 1-convex if $A > 0$ is positive definite, while $\partial M$ is $n$-convex if $H > 0$. 
It is easy to see that $A$ is $(n-1)$-convex if and only if the form $H\gamma - A$ is 
positive definite, 
\begin{equation}\label{e3.23}
H\gamma - A > 0.
\end{equation}
This condition is equivalent to the local convexity of $\partial M$ in $(M, g)$ when 
$n = 2$, but becomes progressively weaker in higher dimensions. 

\begin{proposition}\label{p3.5}
If $\partial M$ is $(n-1)$-convex, so that \eqref{e3.23} holds, then the space 
$${\mathcal V} =  \overline{Im D\Pi_{D}},$$
is of finite codimension in $S_{2}^{m,\alpha}(\partial M)$, where the closure is 
taken in the $C^{m-1,\alpha}$ topology.
\end{proposition}

{\bf Proof:} Recall from Proposition 3.2 that the operator $L$ in \eqref{e2.16} with 
boundary data
\begin{equation}\label{e3.24}
\beta (h) = 0, \ [h^{T}]_{B} = h_{1}, \ tr_{\sigma}A'_{h} = h_{2},
\end{equation}
is elliptic, of Fredholm index 0, provided $\sigma$ is positive (or negative) 
definite and provided $\tau_{B} = B - (tr_{\gamma}B)\gamma \in S_{2}^{m,\alpha}(\partial M)$,
is also negative definite. For $B = A$, by \eqref{e3.23} one has the required 
definiteness, but there is a loss of one derivative in that $\tau \in S_{2}^{m-1,\alpha}
(\partial M)$. Thus, let $A_{\varepsilon}$ be a ($C^{\infty}$) smoothing of $A$, 
$\varepsilon$-close to $A$ in the $C^{m-1,\alpha}$ topology. Then the system $L(h) = 0$ 
with boundary data 
\begin{equation}\label{e3.25}
\beta (h) = 0, \ [h^{T}]_{A_{\varepsilon}} = h_{1}, \ tr_{\sigma}A'_{h} = h_{2},
\end{equation}
is elliptic, of Fredholm index 0. The kernel and cokernel are of finite and equal 
dimensions. 

  Let $\pi_{A_{\varepsilon}}$ denote the projection onto $S_{2}^{m,\alpha}
(\partial M)/\langle A_{\varepsilon}\rangle = S_{\varepsilon}^{m,\alpha}(\partial M)$. 
Then the image of $\pi_{A_{\varepsilon}}\circ D\Pi_{D}$ is of finite codimension in 
$S_{A_{\varepsilon}}^{m,\alpha}(\partial M)$. The fiber 
$(\pi_{A_{\varepsilon}})^{-1}(0)$ consists of symmetric forms of the form 
$fA_{\varepsilon}$. Note that the forms $fA$ are in $Im D\Pi_{D}$, when $\Pi_{D}$ 
is extended to the domain ${\mathbb E}^{m-1,\alpha}(M)$, in that
$$fA = \delta^{*}(fN) \ \ {\rm at} \ \ \partial M,$$
where $\delta^{*}(fN)$ is extended to $M$ to be in Bianchi gauge. Since the forms 
$fA$ are $C\varepsilon$-close to $fA_{\varepsilon}$ in the $C^{m-1,\alpha}$ topology, 
when $|f|_{C^{m-1,\alpha}} \leq C$, it follows (by letting $\varepsilon \rightarrow 
0$) that the closure of $Im D\Pi_{D}$ is of finite codimension in $S_{2}^{m,\alpha}
(\partial M)$. 

{\endproof}

\begin{remark}\label{r3.6}
{\rm Consider hyperbolic 3-space ${\mathbb H}^{3}(-1)$ divided by translation along a 
geodesic, giving a hyperbolic metric $g_{-1}$ on $D^{2}\times S^{1}$. The metric 
$g_{-1}$ has the simple form
$$g_{-1} = dr^{2} + \sinh^{2}r (d\theta_{1})^{2} + \cosh^{2}r (d\theta_{2})^{2}.$$
As in the example discussed in the Introduction, let $\sigma$ be any smooth embedded 
closed curve in the hyperbolic plane $D^{2} = {\mathbb H}^{2}(-1)$ surrounding the 
origin and let $D$ be the disc bounded by $\sigma$. Let 
$M = \pi^{-1}(D^{2}) \simeq D^{2}\times S^{1}$ with $\partial M = \pi^{-1}(\sigma)$, 
so that $M$ is a solid torus with boundary a flat torus $T^{2}$. 

   It is easy to see that $\partial M$ is convex in $M$ whenever $\sigma$ is convex 
in ${\mathbb H}^{2}(-1)$. However the flat torus boundary has two Killing fields, only 
one of which (namely the vertical field tangent to $\theta_{2}$) extends to a Killing 
field on $M$ whenever the geodesic or mean curvature of $\sigma$ in ${\mathbb H}^{2}(-1)$ 
is non-constant. Thus, isometry extension fails, even though $\partial M$ is strictly 
convex - in contrast to the case of rigidity of convex surfaces in ${\mathbb R}^{3}$, 
cf.~\cite{Sp}. 
}
\end{remark}

\section{Isometry Extension and the Divergence Constraint}
\setcounter{equation}{0}

  By Proposition 2.5, the basic issue for the isometry extension property is to 
understand when a Killing field on $(\partial M, \gamma)$ preserves the $2^{\rm nd}$ 
fundamental form $A$ of $\partial M$ in $M$. We begin with the following identity 
on $(\partial M, \gamma)$, which holds on any closed oriented Riemannian manifold. 

\begin{proposition} \label{p4.1}
Let $X$ be a Killing field on $(\partial M, \gamma)$. Suppose $\tau$ is a 
divergence-free symmetric bilinear form on $(\partial M, \gamma)$. Then 
\begin{equation}\label{e4.1}
\int_{\partial M}\langle {\mathcal L}_{X}\tau, h \rangle dV_{\gamma} = 
-2\int_{\partial M}\langle \delta'\tau, X \rangle dV_{\gamma},
\end{equation}
where ${\mathcal L}_{X}$ is the Lie derivative with respect to $X$ and $\delta' = 
\frac{d}{ds}\delta_{\gamma + sh}$ is the variation of the divergence on 
$(\partial M, \gamma)$ in the direction $h \in S_{2}(\partial M)$. 
\end{proposition}

{\bf Proof:} Since the flow of $X$ preserves $\gamma$, one has
\begin{equation}\label{e4.2}
\int_{\partial M}\langle {\mathcal L}_{X}\tau, h \rangle dV_{\gamma} =
-\int_{\partial M}\langle \tau, {\mathcal L}_{X}h \rangle dV_{\gamma}.
\end{equation}
Next, setting $\gamma_{s} = \gamma + sh$, the divergence theorem applied to 
the 1-form $\tau (X)$ on $\partial M$ gives
\begin{equation}\label{e4.3}
0 = \int_{\partial M}\delta_{\gamma_{s}}(\tau(X))dV_{\gamma_{s}} = 
\int_{\partial M}\langle \delta_{\gamma_{s}}\tau, X \rangle dV_{\gamma_{s}} - 
{\tfrac{1}{2}}\int_{\partial M}\langle \tau, {\mathcal L}_{X}\gamma_{s} 
\rangle dV_{\gamma_{s}},
\end{equation}
where the second equality is a simple computation from the definitions; the inner 
products are with respect to $\gamma_{s}$. Taking the derivative with respect to 
$s$ at $s = 0$ and using the facts that $X$ is a Killing field on $\partial M$ 
and $\delta \tau = 0$, it follows that 
$$\int_{\partial M}\langle \delta '\tau, X \rangle dV - {\tfrac{1}{2}}
\int_{\partial M}\langle \tau, {\mathcal L}_{X}h \rangle dV = 0.$$
Combining this with \eqref{e4.2} then gives \eqref{e4.1}. 
{\endproof}

  We now examine the right side of \eqref{e4.1} in connection with the divergence 
constraint \eqref{e3.3}; of course \eqref{e3.3} implies that the form 
$$\tau_{A} \equiv \tau = A - H\gamma,$$
cf.~\eqref{e3.8}, is divergence-free on $\partial M$. 

  We first discuss the general perspective. As discussed in Section 2, one may view the pair 
$(\gamma, A)$ as Cauchy data for the Einstein equations \eqref{e2.3} at $\partial M$. 
The data $(\gamma, A)$ are then formally freely specifiable subject to the 
constraints \eqref{e3.3}-\eqref{e3.4}. Let ${\mathcal T}$ be the space of pairs 
$(\gamma, \tau)$ with $\tau$ divergence-free with respect to $\gamma$; here $\gamma \in 
Met^{m,\alpha}(\partial M)$, $\tau \in S_{2}^{m-1,\alpha}(\partial M)$. One has a 
natural projection onto the first factor 
\begin{equation}\label{e4.4}
\pi: {\mathcal T} \rightarrow Met^{m,\alpha}(\partial M),
\end{equation}
Let also ${\mathcal F} \subset {\mathcal T}$ be the subset of pairs satisfying 
the scalar constraint equation \eqref{e3.4}. When expressed in terms of 
$\tau = A - H\gamma$, \eqref{e3.4} is equivalent to 
$$|\tau|^{2} - \frac{1}{n-1}(tr \tau)^{2} + s_{\gamma} = (n-1)\lambda.$$

  Pairs $(\gamma, \tau) \in {\mathcal F}$ determine formal solutions of the Einstein 
equations near $\partial M$. More precisely, let $(t, x^{i})$ be geodesic boundary 
coordinates for $(M, g)$, so that by the Gauss Lemma, the metric $g$ has the form
\begin{equation}\label{e4.5}
g = dt^{2} + g_{t},
\end{equation}
where $t(x) = dist_{g}(x, \partial M)$ and $g_{t}$ is the induced metric on 
the level set $S(t)$ of $t$. Pulling back by the flow lines of $\nabla t$, $g_{t}$ may 
be viewed as a curve of metrics on $\partial M$, and one may formally expand 
$g_{t}$ in its Taylor series:
\begin{equation}\label{e4.6}
g_{t} \sim \gamma - tA - {\tfrac{1}{2}}t^{2}\dot A + \cdots,
\end{equation}
where $\dot A = \nabla_{N}A = -\nabla_{T}A$, $T = \nabla t = -N$. As noted above, 
the terms $(\gamma, A)$ are freely specifiable, subject to the constraints 
\eqref{e3.3}-\eqref{e3.4}. All the higher order terms in the expansion \eqref{e4.6} 
are then determined by $\gamma$ and $A$. To see this, one first uses the standard 
Riccati equation 
\begin{equation} \label{e4.7}
\nabla_{T}A + A^{2} + R_{T} = 0,
\end{equation}
where $R_{T}(X,Y) = \langle R(X, T)T, Y \rangle$, cf.~\cite{P}. A standard formula gives 
$\nabla_{T}A = {\mathcal L}_{T}A - 2A^{2}$. Also, by the Gauss equation, the curvature 
term $R_{T}$ may be expressed as
$$R_{T} = Ric^{T} - Ric_{int} + HA - A^{2},$$
where $H = tr A$, $Ric_{int}$ is the intrinsic Ricci curvature of $S(t)$ and $Ric^{T}$ 
is the tangential part (tangent to $S(t)$) of the ambient Ricci curvature. 
Substituting in \eqref{e4.7} gives
\begin{equation} \label{e4.8}
\ddot g = -2Ric^{T} + 2Ric_{int} + 4A^{2} - 2HA.
\end{equation}
For Einstein metrics satisfying \eqref{e2.3}, the right side of \eqref{e4.8} involves only 
the first order $t$-derivatives of the metric $g$. Thus, repeated differentiation of 
\eqref{e4.8} shows that all derivatives $g_{(k)} = {\mathcal L}_{T}^{k}g$ are determined 
at the boundary $M$ by the Cauchy data $(\gamma, A)$, so that $(\gamma, A)$ determines 
the formal Taylor expansion of the curve $g_{t}$ in \eqref{e4.5} at $t = 0$. 

  The Cauchy-Kovalevsky theorem implies that if $(\gamma, \tau)$ are real-analytic 
forms on $\partial M$, then the formal series \eqref{e4.6} converges to $g_{t}$, so 
that one obtains an actual Einstein metric $g$ as in \eqref{e4.5}, defined 
in a neighborhood of $\partial M$. Of course, such metrics will not in general 
extend to globally defined Einstein metrics on $M$. 

\medskip

   Now the right side of \eqref{e4.1} is closely related to the linearization of the 
divergence constraint. Thus, if $(\gamma_{s}, \tau_{s})$ is a curve in $Met^{m,\alpha}
(\partial M)\times S_{2}^{m-1,\alpha}(\partial M)$ with tangent vector $(\gamma', \tau') 
= (h, \tau')$ at $s = 0$, then by the Gauss-Codazzi equation one has
\begin{equation}\label{e4.9}
\delta'(\tau) + \delta(\tau') = -(Ric(N, \cdot))' ,
\end{equation}
where $\delta'$ is defined as in \eqref{e4.1}. If $(\gamma_{s}, \tau_{s})$ is a curve 
in ${\mathcal T}$, then
\begin{equation} \label{e4.10} 
\delta'(\tau) + \delta(\tau') = 0;
\end{equation}
this is the linearized divergence constraint. 
\begin{lemma}\label{l4.2}
If the derivative $D\pi$ in \eqref{e4.4} is surjective at $(\gamma, \tau)$, 
$\tau = A - H\gamma$, then
\begin{equation}\label{e4.11}
{\mathcal L}_{X}A = 0 \ \ {\rm on} \ \ \partial M,
\end{equation}
for any Killing field $X$ on $(\partial M, \gamma)$. Conversely, if \eqref{e4.11} 
holds for all such Killing fields $X$, then $D\pi$ is surjective. 
\end{lemma}
{\bf Proof:}
This result follows easily from Proposition 4.1, with $\tau = A - H\gamma$. Thus, 
\eqref{e4.10} gives $\delta'(\tau) = -\delta(\tau ')$, for the variation 
$\delta'$ of $\delta$ in any direction $h \in T_{\gamma}Met(\partial M)$, for 
some $\tau'$. Hence, \eqref{e4.1} gives
\begin{equation}\label{e4.12}
{\mathcal F}(h) = \int_{\partial M}\langle {\mathcal L}_{X}\tau, h \rangle = 
-2\int_{\partial M}\langle \delta(\tau'), X \rangle = 2\int_{\partial M}\langle 
\tau ', \delta^{*}X \rangle = 0,
\end{equation}
since $X$ is a Killing field on $(\partial M, \gamma)$. Since $h$ is arbitrary, 
this implies that
$${\mathcal L}_{X}\tau = 0,$$
on $\partial M$, and \eqref{e4.11} follows by taking the trace of this equation. 
The same proof also gives the converse as well, using the splitting \eqref{e4.13} 
below. 
{\endproof}

   Thus, given $g \in {\mathbb E}$ and its corresponding $2^{\rm nd}$ fundamental 
form $A$, giving the pair $(\gamma, A)$ at $\partial M$, a fundamental issue is whether 
$D\pi$ is surjective at $(\gamma, A)$, i.e.~whether the linearized divergence 
constraint \eqref{e4.10} is solvable, for any variation $h$ of $\gamma$ on 
$\partial M$ (or for a space of variations dense in $S_{2}(\partial M)$ in 
the $L^{2}$ norm). One cannot expect that this holds at a general pair 
$(\gamma, \tau) \in {\mathcal T}$. Namely, for any compact manifold $\partial M$, 
one has 
\begin{equation} \label{e4.13}
\Omega^{1}(\partial M) = Im \delta \oplus Ker \delta^{*}, 
\end{equation}
where $\Omega^{1}$ is the space of ($C^{m-1,\alpha}$) 1-forms on $\partial M$. 
Thus, solvability at $(\gamma, \tau)$ in general requires that 
\begin{equation}\label{e4.14}
\delta'(\tau) \in Im \delta = (Ker \delta^{*})^{\perp}.
\end{equation}
Of course $Ker \delta^{*}$ is exactly the space of Killing fields on $(\partial M, 
\gamma)$, and so this space serves as a potential obstruction space.  

  Obviously, $\pi$ is locally surjective when $(\partial M, \gamma)$ has no Killing 
fields. On the other hand, it is easy to construct examples where $(\partial M, \gamma)$ 
does have Killing fields and $\pi$ is not locally surjective. 

\begin{Example} \label{ex4.3}
{\rm  Let $(\partial M, \gamma)$ be a flat metric on the $n$-torus $T^{n}$; 
for example $\gamma = d\theta_{1}^{2} + \cdots +  d\theta_{n}^{2}$. Let $\tau = 
f(\theta_{1})d\theta_{2}^{2}$ (for example). Then $\delta \tau = 0$, for any $C^{1}$ 
function $f(\theta_{1})$. The pair $(\gamma, \tau)$ is in ${\mathcal T}$, 
and in fact in ${\mathcal F} \subset {\mathcal T}$. Letting $X$ be the 
Killing field $\partial_{\theta_{1}}$, one has ${\mathcal L}_{X}\tau \neq 0$ 
whenever $f$ is non-constant, so that by the converse of Lemma 4.2, $\pi$ 
is not locally surjective at such $(\gamma, \tau)$. 

  If $(\gamma, \tau)$ above are real-analytic, then $(\partial M, \gamma)$ is the 
boundary metric of an Einstein metric defined on a thickening $\partial M \times 
I$ of $\partial M$. Of course in general, such thickenings will not extend to 
Einstein metrics on a compact manifold bounding $\partial M$. 

  To obtain examples on compact manifolds, one may use the examples of 
${\mathbb R}^{2}\times S^{1}$, $S^{3}$ or ${\mathbb H}^{3}(-1)/{\mathbb Z}$ 
discussed in the Introduction and Remark 3.5. Here one has an infinite dimensional 
space of isometric embeddings of a flat torus in ${\mathbb R}^{2}\times S^{1}$, $S^{3}$ 
or ${\mathbb H}^{3}(-1)/{\mathbb Z}$ for which Killing fields on the boundary do not 
extend to Killing fields of the ambient space. }
\end{Example}

   Now clearly $D\pi$ is surjective onto $Im D\Pi_{D}$, since $Im D\Pi_{D}$ consists 
of variations of the boundary metric determined by global variations of the Einstein 
metric $g$ on $M$ which of course satisfy \eqref{e4.10}. Hence if $D\Pi_{D}$ is onto, 
or has dense range in $S_{2}(\partial M)$, then Lemma 4.2 holds, i.e.~\eqref{e4.11} 
holds; compare with Remark 3.3. On the other hand, the examples above show that 
whether \eqref{e4.11} holds or not must depend either on global properties of 
$(M, g)$ or extrinsic properties of $\partial M \subset M$.

\medskip

  Next, we place the discussion above in a broader context of rigidity issues. The 
boundary $(\partial M, \gamma)$ of the Einstein manifold $(M, g)$ is called 
infinitesimally (Einstein) rigid if the kernel $K$ of $D\Pi_{D}$ in \eqref{e3.5} is 
trivial, i.e.~$K = 0$. Thus, infinitesimal rigidity is equivalent to the injectivity 
of $D\Pi_{D}$. It is also equivalent to the local rigidity of $(\partial M, \gamma)$ 
(i.e.~the local uniqueness of an Einstein filling $(M, g)$ up to isometry) by 
the manifold theorem, Theorem 2.1. 

   Suppose $X$ is an infinitesimal isometry at $(\partial M, \gamma)$, in that 
$(\delta^{*}X)^{T} = 0$ at $\partial M$ ($X$ is not necessarily tangent to 
$\partial M$). Then as discussed in Remark 3.1, the deformation $\delta^{*}X$ 
may be extended uniquely to $M$ by choosing it to be in Bianchi gauge. Then 
$\delta^{*}X \in K$ and infinitesimal rigidity of $\partial M$ implies that $k = 0$, 
so that $X$ is a Killing field on $(M, g)$. Rigidity in this more restricted sense 
will be called infinitesimal isometric rigidity. Both forms of such rigidity are 
of course generalizations of the isometry extension property discussed in the 
Introduction. 

  One may obtain analogs of Proposition 4.1 and Lemma 4.2 in this context via the 
Einstein-Hilbert action. Thus, recall that Einstein-Hilbert action with 
Gibbons-Hawking-York boundary term on $M$ is 
\begin{equation} \label{e4.15}
I(g) = I_{EH}(g) = -\int_{M}(s_{g} - 2\Lambda)dV_{g} - 2\int_{\partial M}Hdv_{\gamma}, 
\end{equation}
where $\Lambda  = \frac{n-1}{2}\lambda$, cf.~\cite{Ha}. The $1^{\rm st}$ variation 
of $I$ in the direction $h$ is given by
\begin{equation} \label{e4.16}
\frac{d}{dr}I(g+rh) = \int_{M}\langle \hat E_{g}, h \rangle dV_{g} + \int_{\partial M}
\langle \tau_{g}, h \rangle dv_{\gamma}, 
\end{equation}
where $\hat E$ is the Einstein tensor, 
\begin{equation}\label{e4.17a}
\hat E_{g} = Ric_{g} - \frac{s}{2}g + \Lambda g,
\end{equation}
and $\tau  = A - H\gamma$ is as above. Here and below, all parameter 
derivatives are taken at 0. Einstein metrics with $Ric_{g} - \lambda g = 0$ are 
critical points of $I$, among variations vanishing on $\partial M$. 
Consider a $2$-parameter family of metrics $g_{r,s} = g + rh + sk$ 
where $E_{g} = 0$. Then
\begin{equation} \label{e4.17}
\frac{d^{2}}{dsdr}I(g_{r,s}) = \frac{d^{2}}{drds}I(g_{r,s}). 
\end{equation}
Computing the left side of \eqref{e4.17} by taking the derivative of \eqref{e4.16} 
in the direction $k$ gives
\begin{equation} \label{e4.18}
\frac{d^{2}}{dsdr}I(g_{r,s}) = \int_{M}\langle \hat E'(k), h \rangle dV_{g} + 
\int_{\partial M}\langle \tau'_{k} + a(k^{T}), h^{T}\rangle dv_{\gamma}. 
\end{equation}
Since $\hat E_{g} = 0$, there are no further derivatives of the bulk integral in 
\eqref{e4.16}. Also, $a(k) = -2\tau\circ k + \frac{1}{2}(tr_{\gamma}k)\tau$ 
arises from the variation of the metric and volume form in the direction $k$; 
by definition $(\tau\circ k)(V,W) = \frac{1}{2}\{\langle \tau (V), k(W)\rangle 
+ \langle \tau (W), k(V)\rangle \}$. 

   Similarly, for the right side of \eqref{e4.17} one has 
\begin{equation} \label{e4.19}
\frac{d^{2}}{drds}I(g_{r,s}) = \int_{M}\langle \hat E'(h), k\rangle dV_{g} + 
\int_{\partial M}\langle \tau'_{h} + a(h^{T}), k^{T}\rangle dv_{\gamma}. 
\end{equation}
In particular, suppose $k_{D}$ is an infinitesimal Einstein deformation in 
the kernel $K$ from \eqref{e3.5}, so that $k_{D}|_{\partial M} = k^{T} = 0$. 
If $h \in T{\mathbb E}$ is any infinitesimal Einstein deformation, then 
\eqref{e4.17}-\eqref{e4.19} gives, 
\begin{equation} \label{e4.20}
\int_{\partial M}\langle \tau'_{k_{D}}, h\rangle dv_{\gamma} = \int_{\partial M}
\langle \tau'_{h}, k_{D}\rangle dv_{\gamma} = 0. 
\end{equation}
One thus has
$$I'' (k_{D}, h) = 0, $$
on-shell. Note this computation 
recaptures \eqref{e4.12} when $k_{D} = \delta^{*}X = \frac{1}{2}{\mathcal L}_{X}g$. 

\begin{proposition}\label{p4.4}
If $\pi_{1}(M, \partial M) = 0$ and the linearization $D\Pi_{D}$ has dense range in 
$S_{2}^{m,\alpha}(\partial M)$, then $D\Pi_{D}$ is injective, so that $K = 0$ in 
\eqref{e3.5} and $(\partial M, \gamma)$ is infinitesimally Einstein rigid.
\end{proposition}

{\bf Proof:} The proof is a simple consequence of \eqref{e4.16}-\eqref{e4.18} and 
Corollary 2.4. Thus, suppose $k \in K$ so that $k$ is an infinitesimal Einstein 
deformation with $k^{T} = 0$ at $\partial M$. By \eqref{e4.20}, 
\begin{equation}\label{e4.21}
\int_{\partial M}\langle \tau'_{k}, h \rangle = \int_{\partial M}\langle \tau'_{h}, 
k \rangle = 0,
\end{equation}
for any $h \in Im D\Pi_{D}$. Since $Im D\Pi_{D}$ is dense in $S_{2}^{m,\alpha}(\partial M)$, 
it follows that $(\tau'_{k})^{T} = 0$ on $\partial M$. Taking the trace, it follows that
$$k^{T} = 0 \ \ {\rm and} \ \ (A'_{k})^{T} = 0 \ \ {\rm on} \ \ \partial M.$$
It now follows from Corollary 2.4 that $k = \delta^{*}Z$ with $Z = 0$ on $M$, so that 
$k$ is pure gauge. This means that the equivalence class $[k] = 0$ in $T{\mathcal E}$. 
Alternately, assuming without loss of generality that $k$ is in the Bianchi slice 
$\beta (k) = 0$, it follows again from Corollary 2.4 that $k = 0$, which proves the 
result. 

{\endproof}

  It is an open question whether converse holds, i.e.~if the injectivity of $D\Pi_{D}$ 
implies $D\Pi_{D}$ has dense range. By the discussion in Section 3, $D\Pi_{D}$ is never 
surjective onto $S_{2}^{m,\alpha}(\partial M)$, when $m < \infty$.

\begin{remark}\label{r4.5}
{\rm  There are simple examples of Einstein metrics which are not infinitesimally 
rigid, even when $\partial M$ is convex. Perhaps the simplest example is given by 
the curve of Riemannian Schwarzschild metrics $g_{m}$ on ${\mathbb R}^{2} \times S^{2}$, 
given by 
$$g_{m} = V^{-1}dr^{2} + Vd{\theta}^{2} + r^{2}g_{S^{2}(1)},$$
where $V = V(r) = 1 - \frac{2m}{r}$, $r \geq 2m > 0$. Smoothness at the horizon 
$\{r = 2m\}$ requires that $\theta \in [0,\beta]$ where $\beta = 8\pi m$, so that 
$g_{m}$ may be rewritten in the form
\begin{equation}\label{e4.22}
g_{m} = V^{-1}dr^{2} + 64\pi^{2}m^{2}Vd{\theta}^{2} + r^{2}g_{S^{2}(1)},
\end{equation}
where now $\theta \in [0,1]$. This is a curve of complete Ricci-flat metrics, but the 
metrics $g_{m}$ differ from each other just by rescalings and diffeomorphisms. Taking 
the derivative with respect to $m$ gives an infinitesimal Einstein deformation $\kappa$ 
of $g_{m}$:
\begin{equation}\label{e4.23}
\kappa = 2m64\pi^{2}[1 - \frac{3m}{r}]d\theta^{2} + 
\frac{2}{r}(1 - \frac{2m}{r})^{-2}dr^{2}.
\end{equation}
For the moment, fix $m > 0$ and let $M = M(R) = \{2m \leq r \leq R\}$. The restriction 
of $g_{m}$ to $M$ gives a curve of Einstein metrics on the bounded domain $D^{2}\times 
S^{2}$ with boundary $\partial M \simeq S^{1}\times S^{2}$ and boundary metric
$$\gamma = \gamma_{R} =  64\pi^{2}m^{2}[1 - \frac{2m}{R}]d{\theta}^{2} + R^{2}
g_{S^{2}(1)}.$$
Let $\omega (R)$ be the ratio of the radii of the $S^{1}$ and $S^{2}$ factors at 
$\partial M$, so that
$$\omega (R) = \frac{64\pi^{2}m^{2}[1 - \frac{2m}{R}]}{R^{2}}.$$
Then $\omega (R) \rightarrow 0$ as $R \rightarrow 0$ and $R \rightarrow \infty$, 
and has a single maximum value $64\pi^{2}/27$ at the critical point $R = 3m$ 
where $\kappa^{T} = 0$. At this critical radius, equal to the photon radius 
of the Lorentzian Schwarzschild metric, the boundary metric has the form
$$\gamma = \frac{64}{3}\pi^{2}m^{2}d\theta^{2} + 9m^{2}g_{S^{2}(1)},$$
and a simple calculation shows that the $2^{\rm nd}$ fundamental form $A$ 
is umbilic, with 
$$A = \frac{1}{3\sqrt{3}m}\gamma.$$

  The discussion above shows that the Einstein metric $g_{m}$ is not infinitesimally 
rigid on the domain $M(3m)$; the form $\kappa$ in \eqref{e4.23} is in $Ker D\Pi_{D}$. 
Proposition 4.4 implies that $D\Pi_{D}$ does not have dense range on $M(3m)$; in fact 
boundary metrics for which the mass-independent ratio $\omega > \omega_{0} = 64\pi^{2}/27$ 
are not in $Im \Pi_{D}$ (at least along the Schwarzschild curve). The Dirichlet 
boundary map $\Pi_{D}$ has a simple fold behavior near the critical radius, and 
so has local degree 0. It is shown in \cite{Yo} that the Schwarzschild metric 
$g_{m}$ on $M(R)$ is stable, in that the $2^{\rm nd}$ variation of the action 
\eqref{e4.15} is positive definite, for $R < 3m$, while it becomes unstable, 
has a negative mode or eigenvalue) when $R > 3m$. 

    A detailed discussion of the physical aspects of the Schwarzschild curve is 
given in \cite{Yo}, and further examples in both four and higher dimensions are 
discussed in \cite{AG} and references therein. 

  A simple computation using \eqref{e2.11} shows that on the domain $M(3m)$
$$A'_{\kappa} = \frac{1}{\sqrt{3}m^{2}}(\theta^{1})^{2} - \frac{1}{3\sqrt{3}m^{2}}
\gamma,$$
where $\theta^{1}$ is the unit 1-form in the direction $\theta$. This shows that 
$H'_{\kappa} = 0$ at $\partial M$. Hence, the form $\kappa$ is also in the kernel 
of the Fredholm boundary map $\widetilde \Pi_{B,\gamma}$ in \eqref{e3.17}. This 
shows that the generalization of Theorem 1.1 to infinitesimal Einstein rigidity 
is false; the form $\kappa$ is a non-trival infinitesimal Einstein deformation 
preserving the boundary metric and mean curvature. Of course $\kappa$ is not of 
the form $\delta^{*}X$ for some vector field $X$. }
\end{remark}

\begin{remark}\label{r4.6}
{\rm  The work above has been carried out with the operator $E(g) = 
Ric_{g} - \lambda g$ in \eqref{e2.7}, and in the Bianchi gauge, since the 
computations are the simplest in this setting. However, the discussion 
following \eqref{e4.16} suggests that the ``physical'' Einstein operator 
$$\hat E(g) = Ric_{g} - \frac{s}{2}g + \Lambda g$$
in \eqref{e4.17a} may be more natural in certain respects. This is in fact 
the case, and will be used in Section 5. 

  To set the stage, note first that the analog of the Bianchi identity 
in this setting is $\delta \hat E = 0$. As in Section 2, fix any background 
Einstein metric $\widetilde g$ and consider the operator
\begin{equation}\label{e4.24}
\hat \Phi_{\widetilde g}(g) = Ric_{g} - \frac{s}{2}g + \Lambda g + 
\delta_{g}^{*}\delta_{\widetilde g}(g),
\end{equation}
(cf.~\eqref{e2.15}). The linearization of $\hat \Phi$ at $g = \widetilde g$ is 
\begin{equation}\label{e4.25}
\hat L(h) = \nabla^{*}\nabla h - 2R(h) - (D^{2}tr h + \delta \delta h)g + \Delta tr h \,g,
\end{equation}
where $D^{2}$ is the Hessian and $\Delta = tr D^{2}$ the Laplacian (with respect to $g$). 
This is more complicated than \eqref{e2.16}, but it is easy to see that $\hat L$ 
is formally self-adjoint; this also follows directly from the symmetry of the 
$2^{\rm nd}$ derivatives in \eqref{e4.17}-\eqref{e4.18}. Of course solutions 
of $\hat L(h) = 0$ with $\delta(h) = 0$ on $M$ are infinitesimal Einstein 
deformations. 

  It is straightforward to see that the operator $\hat L$ with boundary conditions
\begin{equation}\label{e4.26}
\delta(h) = 0, \ [h^{T}]_{0} = h_{1}, \ H'_{h} = h_{2} \ \ {\rm at} \ \ \partial M,
\end{equation}
(where $[h^{T}]_{0} = [h^{T}]{\gamma}$, the trace-free part of $h^{T}$) is a 
well-posed elliptic boundary value problem. This follows from Proposition 3.2, 
using the fact that the change from Bianchi to divergence-free gauge is unique; 
it is also proved directly in \cite{AK}. Moreover, it is not only of Fredholm 
index 0, but the boundary value problem \eqref{e4.26} is self-adjoint (which 
is not the case for $L$ with the boundary conditions \eqref{e3.10} from 
Proposition 3.2) cf.~again \cite{AK}.  

  All of the remaining discussion in Section 2 carries over immediately to this setting; 
the only change is that one replaces the Bianchi operator $\beta$ by the divergence 
operator $\delta$. For instance, the analog of \eqref{e2.17} is
\begin{equation}\label{e4.27}
L_{\hat E} = \hat L - 2\delta^{*}\delta.
\end{equation}
Similarly, Lemma 2.3 holds in this setting, with the same proof. 

}
\end{remark}

\section{Proof of the Main Results.}

\setcounter{equation}{0}

  In this section, we prove the main results discussed in the Introduction, beginning 
with Theorem 1.1. As noted above, one needs to use global arguments to prove Theorem 1.1. 
We do this by studying global properties of the linearized operator $\hat L$ from 
\eqref{e4.25}. 

  Consider then the elliptic boundary value problem \eqref{e4.25}-\eqref{e4.26}: 
\begin{equation}\label{e5.1}
\hat L(h) = \ell, \ \ {\rm on} \ M, \ \ \delta(h) = h_{0}, \ [h^{T}]_{0} = h_{1}, 
\ H'_{h} = h_{2} \ {\rm on} \ \partial M,
\end{equation}
As noted in Remark 4.6, this is a self-adjoint elliptic boundary value problem. The 
self-adjoint property leads to significant simplifications in the proof, which is 
why we use the divergence gauge and $\hat L$ in place of $L$ and the Bianchi gauge. 

   Let $K$ denote the kernel, so that $k \in K$ means 
$$\hat L(k) = 0, \ \ \delta(k) = 0, \ [k^{T}]_{0} = 0, \ H'_{k} = 0.$$
If $K = 0$, then $L$ is surjective and so the form $\ell$ and boundary values for 
$\delta(h), [h^{T}]_{0}$ and $H'_{h}$ may be freely chosen; given arbitrary $\ell$ 
and $h_{i}$, $0 \leq i \leq 2$, the system \eqref{e5.1} has a unique solution, 
(when suitable smoothness assumptions are imposed).  

  Now (regardless of whether $K = 0$ or not) as in Lemma 4.2 (using \eqref{e4.10} 
and \eqref{e4.12}) one has
\begin{equation}\label{e5.2}
\int_{\partial M}\langle {\mathcal L}_{X}\tau, h^{T} \rangle = 
\int_{\partial M}(Ric(N, X))'_{h}.
\end{equation}
We will prove that any deformation $h^{T}$ of $\gamma$ on $\partial M$ extends 
to a deformation $h$ of $g$ on $M$ such that the right side of \eqref{e5.2} 
vanishes; Theorem 1.1 then follows easily via Proposition 2.5.  

  Note first that \eqref{e5.2} vanishes in pure-trace directions $h^{T} = 
f \gamma$. Namely, since $X$ is Killing, $tr ({\mathcal L}_{X}\tau) = 
-(n-1)X(H) = 0$, by assumption. Hence, $\langle {\mathcal L}_{X}\tau, f \gamma 
 \rangle = 0$ pointwise and so the right side of \eqref{e5.2} vanishes in 
pure-trace directions also.  

  By Lemma 2.3 and Remark 4.6, deformations $h$ satisfying 
\begin{equation}\label{e5.3}
\hat L(h) = 0, \ \ {\rm on} \ M, \ \ \delta(h) = 0 \ \ {\rm on} \ \partial M,
\end{equation}
are infinitesimal Einstein deformations in divergence-free gauge on $M$ and hence, 
at $\partial M$,
$$(Ric(N, X))'_{h} = 0,$$
since $N$ is normal and $X$ is tangential. Now write any $h^{T}$ on $\partial M$ as 
$h^{T} = h_{0} + f \gamma$ where $h_{0}$ is trace-free. Let $\bar f$ be any 
smooth function and let $\bar h^{T} = h_{0} + \bar f \gamma$, so that 
$\bar h^{T} - h^{T}$ is pure-trace. Then by the remarks following \eqref{e5.2}, 
\begin{equation}\label{e5.4}
\int_{\partial M}(Ric(N, X))'_{\bar h} = 
\int_{\partial M}(Ric(N, X))'_{h}.
\end{equation}

  Suppose first the boundary value problem in \eqref{e5.1} has trivial kernel, 
$K = 0$. It follows that there exists an infinitesimal Einstein deformation $h$ 
of $(M, g)$ satisfying \eqref{e5.3} with $h^{T} = h_{0} + f \gamma$, for some 
$f$ and with the class $[h^{T}]_{0} = h_{1}$ arbitrarily prescribed. For all 
such $h$, it follows that
\begin{equation}\label{e5.5}
\int_{\partial M}(Ric(N, X))'_{h} = 0.
\end{equation}
Via \eqref{e5.4}, \eqref{e5.5} then also holds for all $h$, and so by \eqref{e5.2}, 
one obtains
$${\mathcal L}_{X}\tau = 0.$$
Since $tr({\mathcal L}_{X}\tau) = 0$, this gives ${\mathcal L}_{X}A = 0$ and 
Theorem 1.1 then follows from Proposition 2.5.

  Next, suppose $K \neq 0$. Let $S_{0}^{m,\alpha}(M)$ be the Banach space of symmetric 
forms on $M$ with 0 boundary values in \eqref{e5.1}. Let $K^{\perp}$ be the $L^{2}$ 
orthogonal complement of $K$ within $S_{0}^{m,\alpha}(M)$. This is a closed subspace 
of $S_{0}^{m,\alpha}(M)$, of finite codimension with complement $K$, so that 
$S_{0}^{m,\alpha}(M) = K^{\perp}\oplus K$. The operator $\hat L|_{K^{\perp}}$ is an 
isomorphism onto its image $Im(\hat L)$, and since $\hat L$ is self-adjoint, 
$Im(\hat L) = K^{\perp}$. The kernel $K$ is the orthogonal slice to the image, 
and 
\begin{equation}\label{e5.6}
\delta (k) = 0,
\end{equation}
for all $k \in K$. 

  We now construct a different linear slice $\widetilde Q$ to $Im(\hat L)$ with certain 
specific properties at $\partial M$, which are not known to hold, apriori, for $K$. 
First choose a basis $\{k_{i}\}$ for $K$ and choose $C^{m,\alpha}$ symmetric forms $q_{i}$ 
of compact support in $M$ such that, for each $i$,
\begin{equation}\label{e5.7}
\int_{M}\langle q_{i}, k_{i} \rangle \neq 0.
\end{equation}
The span of $\{q_{i}\}$ gives a linear space $Q \simeq K$, with $Q$ nowhere orthogonal 
to $K$, i.e.~no form $q \in Q$ is orthogonal to $K$, so that $Q$ is also a slice to 
$Im(\hat L)$. Next, for a fixed $q = q_{i}$, consider the forms $q + \delta^{*}Y$, where 
$Y$ is a solution of the equation
\begin{equation}\label{e5.8}
\delta \delta^{*}Y = -\delta(q),
\end{equation}
so that $\delta(q + \delta^{*}Y) = 0$. Now one can solve the equation \eqref{e5.8} with 
either Dirichlet, Neumann or mixed (Robin) boundary conditions at $\partial M$. The 
two boundary conditions we impose are:
\begin{equation}\label{e5.9}
\int_{\partial M}\langle k(N), Y\rangle = 0,
\end{equation}
\begin{equation}\label{e5.10}
\int_{\partial M}(\delta^{*}Y)(N,X) = \int_{\partial M}\langle \nabla_{N}Y - A(Y), 
X \rangle = 0.
\end{equation}
Here $X$ is the given Killing field on $\partial M$. The first equality in \eqref{e5.10} 
follows directly from the definition of $\delta^{*}Y$ (using the fact that 
$\int_{\partial M}X(f) = \int_{\partial M}f\delta X = 0$, since $X$ is Killing on 
$\partial M$) so only the second equality is a condition. 

  These are mixed Dirichlet-Neumann conditions, which are straightforward to solve. 
In detail, consider first the homogeneous equation
\begin{equation}\label{e5.10a}
\delta \delta^{*}Y = 0.
\end{equation}
Let $R_{0}$ be the Dirichlet-to-Robin type map sending Dirichlet data $Y \in 
\chi^{m+1,\alpha}(\partial M)$ (the space of vector fields at $\partial M$) 
to $\nabla_{N}Y - A(Y)$ on $\partial M$, where $Y$ solves \eqref{e5.10a} on $M$; 
$$R_{0}(Y) = \nabla_{N}Y - A(Y).$$ 
The map $R_{0}$ is Fredholm, of Fredholm index 0, with kernel ${\mathcal K}$ equal 
to the space of vector fields $Y$ at $\partial M$ which extend to Killing fields on 
$(M, g)$. This follows by pairing \eqref{e5.10a} with $Y$ and integrating by parts. 
Orthogonal to the kernel, $R_{0}$ is an isomorphism onto its image ${\mathcal V}_{0} 
\subset \chi^{m,\alpha}(\partial M)$ and ${\mathcal V}_{0}\oplus {\mathcal K} = 
\chi^{m,\alpha}(\partial M)$. The Dirichlet-to-Robin map $R_{q}$ for \eqref{e5.8} is 
then an affine map onto the image 
$${\mathcal V}_{q} = {\mathcal V}_{0} + z_{q},$$ 
where $z_{q} = \nabla_{N}Y$ and $Y$ solves \eqref{e5.8} with zero Dirichlet 
boundary data. 

  On the other hand, the condition \eqref{e5.9} defines a codimension 1 hypersurface 
${\mathcal S} \subset \chi^{m+1,\alpha}(\partial M)$ (with normal vector $k(N)$) 
which maps under $R_{q}$ to a codimension 1 hypersurface $R_{q}({\mathcal S})$ of 
${\mathcal V}_{q}$. 

  Suppose first $z_{q} \in {\mathcal V}_{0}$ (e.g.~$(M, g)$ has no Killing fields) 
so ${\mathcal V}_{q} = {\mathcal V}_{0}$. We have then two codimension 1 hypersurfaces 
of ${\mathcal V}_{q}$, namely $R_{q}({\mathcal S})$ and the hypersurface 
${\mathcal T}_{X}$ defined by \eqref{e5.10} with normal vector $X$. Any vector 
field $Y$ such that $R_{q}(Y)$ lies in the intersection of these two hypersurfaces 
satisfies \eqref{e5.9}-\eqref{e5.10}. Since this intersection is of codimension 2 
in ${\mathcal V}_{0}$, it is clear there is a large space of solutions. 

  If however $z_{q} \notin {\mathcal V}_{0}$, then ${\mathcal V}_{q}$ is an 
affine subspace, of finite codimension in $\chi^{m,\alpha}(\partial M)$ with 
$R_{q}({\mathcal S})$ of codimension 1 in ${\mathcal V}_{q}$. Let $z_{q}'$ be 
the vector normal to ${\mathcal V}_{0}$ such that ${\mathcal V}_{q} = {\mathcal V}_{0} 
+ z_{q}'$. Then \eqref{e5.9}-\eqref{e5.10} has no solutions, i.e.~$R_{q}({\mathcal S}) 
\cap {\mathcal T}_{X} = \emptyset$, if and only if the normal vector $z_{q}'$ is a 
constant multiple of the normal vector $X$, so that the functional 
$\int \langle \cdot, X \rangle$ is constant on $R_{q}({\mathcal S})$. 
However, if $X \in {\mathcal K}$, then Theorem 1.1 is proved, and so, without 
loss of generality, one may assume $X \perp {\mathcal K}$, so that $X \in 
{\mathcal V}_{0}$. The functional $\int \langle \cdot, X \rangle$ is then 
non-trivial (i.e.~non-constant) and assumes the value 0 again on a 
codimension 2 subspace of ${\mathcal V}_{q}$. In this way, we see that 
\eqref{e5.9} and \eqref{e5.10} always have a large space of solutions. 

  We pick such a solution $Y_{i}$, for each $q_{i}$ in a basis of $Q$ and, extending linearly, 
let $\widetilde Q = \{\widetilde q = q + \delta^{*}Y\}$, so that $\widetilde Q \simeq K$. 
Observe that $\widetilde Q$ is still not orthogonal to $K$, i.e.~for any $\widetilde q$ 
there exists $k \in K$ such that
\begin{equation}\label{e5.11}
\int_{M}\langle \widetilde q, k \rangle \neq 0.
\end{equation}
It suffices to verify \eqref{e5.11} on a basis $\widetilde q_{i}$ and via \eqref{e5.7} it 
then suffices to show that $\int_{M}\langle \delta^{*}Y, k \rangle = 0$. This follows 
from a standard integration-by-parts:
$$\int_{M}\langle \delta^{*}Y, k \rangle = \int_{M}\langle Y, \delta(k) \rangle  
+ \int_{\partial M}\langle k(N), Y\rangle = 0,$$
where we have used \eqref{e5.6} and \eqref{e5.9}. 

  Thus, $\widetilde Q$ is also not orthogonal to $K$, so gives a slice to 
$Im(\hat L)$, and by \eqref{e5.8},
\begin{equation}\label{e5.12}
\delta(\widetilde q) = 0,
\end{equation}
on $M$, for each $\widetilde q \in \widetilde Q$. 

   Next form the operator 
\begin{equation}\label{e5.13}
\widetilde L(h) = \hat L(h) + \pi_{\widetilde Q}(h),
\end{equation}
where $\pi_{\widetilde Q}$ is the orthogonal projection onto $\widetilde Q$. 
Since by construction $\widetilde Q$ is linearly independent from $Im(\hat L)$, 
it follows that $\widetilde L$ is an isomorphism
\begin{equation}\label{e5.14}
\widetilde L:S_{0}^{m,\alpha}(M) \rightarrow S^{m-2,\alpha}(M). 
\end{equation}

   Consider now the boundary value problem
\begin{equation}\label{e5.15}
\widetilde L(h) = 0, \ \  \delta(h) = 0, \ [h^{T}]_{0} = h_{1}, \ H'_{h} = h_{2}.
\end{equation}
Since $\widetilde L$ in \eqref{e5.14} is a bijection, it follows from a standard 
substraction procedure that \eqref{e5.15} has a unique solution, for arbitrary $h_{1}$ 
and $h_{2}$. Namely, take any symmetric form $v$ satisfying the boundary conditions 
in \eqref{e5.15} and extend $v$ to a smooth form on $M$ (arbitrarily but smoothly) 
so that $\widetilde L(v) = w$, for some $w$. Let $h_{0}$ be the unique solution of 
$\widetilde L(h_{0}) = w$ with 0-boundary values, as in \eqref{e5.14}. Then 
$h = v - h_{0}$ solves \eqref{e5.15}. Of course solutions of \eqref{e5.15} are 
not infinitesimal Einstein deformations in general. 

  Next, we claim that for any $h_{1}$ and $h_{2}$, the solution $h$ of \eqref{e5.15} 
satisfies
\begin{equation}\label{e5.16}
\delta (h) = 0,
\end{equation}
on $M$. To prove this, one has $\hat L(h) = \widetilde L(h) - \pi_{\widetilde Q}(h) 
= -\pi_{\widetilde Q}(h)$. By \eqref{e4.27}, one has $\delta \hat L(h) = 2\delta \delta^{*}
(\delta (h))$ (since $\delta \hat L = 0$ by the Bianchi identity) which gives
$$2\delta \delta^{*}(\delta (h)) = -\delta(\pi_{\widetilde Q}(h)).$$
But $\pi_{\widetilde Q}(h) = \widetilde q$ for some $\widetilde q$ and 
$\delta(\widetilde q) = 0$, by \eqref{e5.12}. So
\begin{equation}\label{e5.17}
\delta \delta^{*}(\delta(h)) = 0,
\end{equation}
on $M$. By assumption (in \eqref{e5.15}) $\delta(h) = 0$ on $\partial M$ and so by 
Lemma 2.3 (i.e.~the analog of this result as mentioned in Remark 4.6) 
\eqref{e5.16} follows. 

   The proof of Theorem 1.1 is now easily completed as follows. Namely for $h$ 
as in \eqref{e5.15}, $\delta(h) = 0$ on $M$ and $\hat L(h) = \widetilde L(h) - 
\pi_{\widetilde Q}(h) = -\pi_{\widetilde Q}(h) = -\widetilde q$, for some $\widetilde q$. 
We have $\hat L(h) = \hat E'(h) + 2\delta^{*}\delta(h) = \hat E'(h)$. Thus
\begin{equation}\label{e5.17a}
\hat E'_{h} + \widetilde q = 0.
\end{equation}
Also $\hat E'_{h} = Ric'_{h} - \frac{s'}{2}g - \frac{s}{2}h + \Lambda h$, so that 
$\hat E'_{h}(N,X) = Ric'_{h}(N,X) - \frac{s}{2}h(N,X) + \Lambda h(N,X)$ (since $X$ 
is tangent to $\partial M$). But $Ric'_{h}(N,X) = (Ric(N,X))'_{h} - Ric(N',X) = 
(Ric(N,X))'_{h} - (\frac{s}{2} - \Lambda)\langle N', X \rangle$. Since 
$\langle N', X \rangle = -h(N, X)$, it follows that $\hat E'_{h}(N,X) = 
(Ric(N, X))'_{h}$ on $\partial M$. Recall that $\widetilde q = q + \delta^{*}Y$ 
where $q$ has compact support (so vanishes near $\partial M$) while $\delta^{*}Y$ 
satisfies \eqref{e5.10}, i.e.
\begin{equation}\label{e5.17b}
\int_{\partial M}(\delta^{*}Y)(X,N) = \int_{\partial M}\langle \nabla_{N}Y - 
A(Y), X \rangle = 0.
\end{equation}
This gives
\begin{equation}\label{e5.18}
\int_{\partial M}(Ric(N, X))'_{h} = -\int_{\partial M}\widetilde q(X,N) = 0.
\end{equation}
This holds for $[h^{T}]_{0}$ arbitrary, and so as in \eqref{e5.4}, \eqref{e5.18} holds 
for all variations $h$ on $\partial M$. As above via \eqref{e5.2} and Proposition 2.5, 
this completes the proof.

{\endproof}

\begin{remark}\label{r5.1}
{\rm The proof of Theorem 1.1 above shows that, when for instance $H = const$ at 
$\partial M$, one has
\begin{equation}\label{e5.20}
K \cap Im \delta^{*} = 0,
\end{equation}
where $\delta^{*}$ acts on vector fields $X$ tangent to $\partial M$, and 
$K = Ker D\Pi_{D}$, as in \eqref{e3.5}. 

  However, for instance in dimension 3, all Einstein deformations are pure gauge, 
i.e.~of the form $\delta^{*}V$, for some vector field $V$, not necessarily tangent 
to $\partial M$, cf.~Remark 3.1. Hence, if $\Pi_{D}$ is degenerate at some constant 
curvature metric $(M^{3}, g)$, i.e.~$K = K_{g} \neq 0$ and again $H = const$ at 
$\partial M$, then 
$$K \cap \delta^{*}V \neq 0,$$
for general $V$ at $\partial M$. The condition $H = const$ is necessary here, 
cf.~Example 4.3. }
\end{remark}

  The proof of Theorem 1.1 above generalizes to the case of isometric rigidity, 
where the vector field $X$ is not assumed tangent to $\partial M$, but is a 
general vector field at $\partial M$, preserving the mean curvature. 

\begin{proposition}\label{p5.2}
Let $X$ be a vector field at $\partial M$ generating an infinitesimal isometry at 
$\partial M$ and suppose ${\mathcal L}_{X}(H) = X(H) = 0$. If $\pi_{1}(M, \partial M) 
= 0$, then $X$ extends to a Killing field on $(M, g)$. 
\end{proposition}

{\bf Proof:} The proof is a simple extension of the proof of Theorem 1.1. First, 
it is easy to see that \eqref{e5.17a} holds in general, without the restriction 
that $X$ is tangent to $\partial M$, with a slight redefinition of $\widetilde Q$. 
Namely, for general $X$, the left side of \eqref{e5.10} remains the same, but the 
middle expression has new terms coming from the normal part of $X$; these are 
of the same form as before, and their presence does not affect the validity of 
the proof that \eqref{e5.9}-\eqref{e5.10} are easily solvable for general $X$. 
Next, using \eqref{e4.18}-\eqref{e4.19} together with the fact that 
$(\delta^{*}X)^{T} = 0$, one obtains, for any deformation $h$ of $g$,
\begin{equation}\label{e5.21}
\int_{\partial M}\langle \tau'_{\delta^{*}X}, h \rangle = 
\int_{\partial M}\langle \tau'_{h}, (\delta^{*}X)^{T} \rangle + 
\int_{M}\langle \hat E'_{h}, \delta^{*}X\rangle = 
\int_{M}\langle \hat E'_{h}, \delta^{*}X\rangle.
\end{equation}
Integrating the right-hand side by parts, and using the fact that 
$\delta \hat E' = 0$, it follows that
$$\int_{\partial M}\langle \tau'_{\delta^{*}X}, h \rangle = 
\int_{\partial M}(\hat E'_{h})(X, N).$$
Now as in \eqref{e5.17a}, $(\hat E'_{h})(X, N) = -\widetilde q$, and the analog 
of \eqref{e5.17b} holds as before. Hence, 
\begin{equation}\label{e5.22}
\int_{\partial M}\langle \tau'_{\delta^{*}X}, h \rangle = 0.
\end{equation}

  Since the assumptions ${\mathcal L}_{X}H = 0$ and $(\delta^{*}X)^{T} = 0$ imply 
that $tr \tau'_{\delta^{*}X} = 0$, and since $h$ is arbitrary at $\partial M$ modulo 
pure-trace terms, it follows from \eqref{e5.22} as before in the proof of Theorem 1.1 
that 
$$(\delta^{*}X)^{T} = 0, \ \ {\rm and} \ \ A'_{\delta^{*}X} = 0.$$
at $\partial M$. The result then follows again from Proposition 2.5. 
{\endproof}

{\bf Proof of Corollary 1.2.}

  Theorem 1.1 implies that the isometry group $SO(n+1)$ of $S^{n}$ extends to a group 
of isometries of the Einstein manifold $(M^{n+1}, g)$. This reduces the Einstein 
equations to a simple system of ODE's (the metric $g$ is of cohomogeneity 1) and it 
is standard that the only smooth solutions are given by constant curvature metrics, 
cf.~\cite{Be} for example. 
{\endproof}

  The same proof shows that if $(\partial M, \gamma)$ is homogeneous, then any 
Einstein filling metric $(M, g)$ is of cohomogeneity 1. Such metrics have been 
completely classfied in many situations, cf.~\cite{Be} for further information.

\medskip

  We complete this section with a brief discussion of exterior and global boundary value 
problems. Thus, suppose $M^{n+1}$ is an open manifold with compact ``inner'' boundary 
$\partial M$ and with a finite number of ends, each (locally) asymptotically flat. 
Topologically, each end is of the form $({\mathbb R}^{k} \setminus B)\times T^{n+1-k}$, 
or a quotient of this space by a finite group of isometries. Here $T^{n+1-k}$ is the 
$(n+1-k)$-torus, and we assume $3 \leq k \leq n+1$. Assume also, as usual, that 
$\pi_{1}(M, \partial M) = 0$. An Einstein metric is asymptotically locally flat 
(ALF) if it decays to a flat metric on each end at a rate $r^{-(k-2)}$ (the 
decay rate of the Green's function for the Laplacian) where $r$ is the 
distance from a fixed point. 

  It is proved in \cite{An1} that the analog of Theorem 2.1 holds, namely the space of 
asymptotically locally flat Einstein metrics on an exterior domain $M$ is a smooth 
Banach manifold, for which the Dirichlet boundary map is $C^{\infty}$ smooth. 
Lemma 2.3 also holds in this context. All of the remaining results in Section 2 - 
Section 5 above concern issues at or near $\partial M$, and it is straightforward to verify 
that their proofs carry over to this exterior context without change. In particular 
the analog of Theorem 1.1 holds:

\begin{proposition}\label{p5.3}
Let $g$ be a $C^{m,\alpha}$ Ricci-flat metric on an exterior domain $M$, $m \geq 5$, 
with a finite number of locally asymptotically flat ends. Suppose also \eqref{e1.2} 
holds. Then any Killing field $X$ on $(\partial M, \gamma)$ for which $X(H) = 0$, 
extends uniquely to a Killing field on $(M, g)$. 
\end{proposition}
{\endproof}

  Next we point out that an analog of Theorem 1.1 holds for complete conformally 
compact Einstein metrics, where the boundary is at infinity (conformal infinity). 
The proof below corrects an error in the proof of this result in \cite{AH}. 

\begin{theorem}\label{t5.4}
Let $(M, g)$ be a conformally compact Einstein metric, with smooth conformal 
infinity $(\partial M, [\gamma])$ and suppose $\pi_{1}(M, \partial M) = 0$. 
Then any (conformal) Killing field of $\partial M$ extends to a Killing field 
of $(M, g)$. 
\end{theorem}

{\bf Proof:} The proof is a simple adaptation of the proof of Theorem 1.1, using 
information provided in \cite{AH}, to which we refer for some details. Let 
$t$ be a geodesic compactification of $(M, g)$ and let $S(t)$ and $B(t)$ be the 
level and super-level sets of $t$, so that $\partial M = S(0)$, $M = B(0)$. One 
has
$$g = t^{-2}\bar g = t^{-2}(dt^{2} + \bar g_{t}) = d(\log t)^{2} + g_{t},$$
where $\bar g$ is the conformal compactification of $g$ with respect to $t$. 
The Killing vector field $X$ on $(\partial M, \gamma)$ is extended into $M$ to be 
in Bianchi gauge. As shown in \cite{AH}, $\delta^{*}X$ is then transverse-traceless 
and one has the following estimates: $\langle X, N \rangle = O(t^{n+1})$, 
$[X, N] = O(t^{n+1})$, where $N = -t\partial_{t}$ is the unit outward normal at $S(t)$, 
and $\delta^{*}X = O(t^{n})$, $X(H) = O(t^{n+1})$. 

  Now as in \eqref{e5.21}-\eqref{e5.22} (setting $M = B(t)$) one then has
\begin{equation}\label{e5.24}
{\tfrac{1}{2}}\int_{S(t)}\langle {\mathcal L}_{X}\tau, h \rangle = 
\int_{S(t)}\langle \tau'_{h}, (\delta^{*}X)^{T}\rangle + 
\int_{B(t)}\langle \hat E'_{h}, \delta^{*}X \rangle = \int_{S(t)}\langle \tau'_{h}, 
(\delta^{*}X)^{T}\rangle. 
\end{equation}
The equation \eqref{e5.24} holds for any $h$ such that $\widetilde L(h) = 0$ 
satisfying the boundary conditions $\delta(h) = 0$, $[h^{T}]_{0} = h_{1}$, 
$H'_{h} = h_{2}$ with $h_{1}$, $h_{2}$ arbitrary on $S(t)$. This follows from 
the discussion above concerning \eqref{e5.15}, which holds without the assumption 
that $X$ is Killing on the boundary $S(t)$. As in \cite{AH}, let $\kappa = \delta^{*}X$ and 
$\widetilde \kappa = t^{-n}\delta^{*}X$, so that $\widetilde \kappa$ is uniformly bounded as 
$t \rightarrow 0$; one has $tr \widetilde \kappa = 0$, $\widetilde \kappa(N, \cdot) = O(t)$, 
and $\widetilde \kappa (N, N) = o(t)$. 

  We then choose $h$ such that $[h^{T}]_{0} = [\widetilde \kappa^{T}]_{0}$ 
on $S(t)$. Substituting this in \eqref{e5.24} gives the basic relation
\begin{equation}\label{e5.25}
\int_{S(t)}\langle {\mathcal L}_{X}\tau, \widetilde \kappa \rangle = 
2\int_{S(t)}\langle \tau'_{\widetilde \kappa}, (\delta^{*}X)^{T}\rangle. 
\end{equation}

  On the other hand, we compute both sides of \eqref{e5.25} directly, as in \cite{AH}. 
To set the stage for this, write
$$g_{s} = g + s\kappa + O(s^{2}) = g + s\delta^{*}X + O(s^{2}).$$
If $t_{s}$ is the geodesic defining function for $g_{s}$, (with boundary metric 
$\gamma$), then the Fefferman-Graham expansion gives $\bar g_{s} = dt_{s}^{2} + 
(\gamma + t_{s}^{2}g_{(2),s} + \dots + t_{s}^{n}\log t_{s}{\mathcal H}_{s} + 
t_{s}^{n}g_{(n),s}) + O(t^{n+1})$. One has $t_{s} = t + sO(t^{n+2}) + O(s^{2})$. 
Taking the derivative with respect to $s$ at $s = 0$, and using the fact that $X$ 
is Killing on $(\partial M, \gamma)$, together with the fact that the lower order 
terms $g_{(k)}$, $k < n$, and ${\mathcal H}$ are determined by $\gamma$, 
it follows that on $\partial M$,
\begin{equation}\label{e5.26}
\hat \kappa = {\tfrac{1}{2}}{\mathcal L}_{X}g_{(n)},
\end{equation}
where $\hat \kappa = \lim_{t\rightarrow 0}t^{-(n-2)}\kappa$; here both $\hat \kappa$ 
and ${\mathcal L}_{X}g_{(n)}$ are viewed as forms on $(\partial M, \gamma)$. 

  Now we first claim that as $t \rightarrow 0$, 
\begin{equation}\label{e5.27}
\int_{S(t)}\langle {\mathcal L}_{X}\tau, \widetilde \kappa \rangle = 
-\frac{3n-2}{4}\int_{\partial M}|{\mathcal L}_{X}g_{(n)}|^{2}dV_{\gamma} 
+ o(1).
\end{equation}
To start, on $(S(t), g_{t})$ one has 
\begin{equation}\label{e5.28}
{\mathcal L}_{X}A = -{\tfrac{n-2}{2}}t^{n-2}{\mathcal L}_{X}g_{(n)} +
 O(t^{n-1}),
\end{equation}
To see this, $A = \frac{1}{2}{\mathcal L}_{N}g = -\frac{1}{2}{\mathcal L}_{t\partial t}g 
= -\frac{1}{2}{\mathcal L}_{t\partial t}(t^{-2}\bar g_{t})$. But 
${\mathcal L}_{t\partial t}(t^{-2}\bar g_{t}) = \sum {\mathcal L}_{t\partial t}
(t^{-2+k}g_{(k)}) = \sum (k-2)t^{k-2}g_{(k)}$. Since ${\mathcal L}_{X}g_{(k)} = 0$ 
and ${\mathcal L}_{X}{\mathcal H} = 0$, , it follows that
$$\int_{S(t)}\langle {\mathcal L}_{X}A , \widetilde \kappa \rangle_{g_{t}} dV_{S(t)} = 
-{\tfrac{n-2}{2}}\int_{S(t)}\langle {\mathcal L}_{X}g_{(n)}, \hat \kappa \rangle_{\gamma} 
dV_{\gamma} + O(t).$$
Next, one has ${\mathcal L}_{X}(Hg_{t}) = X(H)g_{t} + H{\mathcal L}_{X}g_{t}$. For the 
first term, $X(H) = tr {\mathcal L}_{X}A + O(t^{n}) = -\frac{n-2}{2}t^{n-2}tr {\mathcal L}_{X}g_{(n)} 
+ O(t^{n})$. Since $tr g_{(n)}$ is intrinsic to $\gamma$ and $X$ is Killing on $(\partial M, \gamma)$, 
it follows that $X(H) = O(t^{n-1})$. Also, $\langle g_{t}, \widetilde \kappa \rangle = -\widetilde 
\kappa(N, N) = o(t)$. Hence $X(H)\langle g_{t}, \widetilde \kappa \rangle dV_{S(t)} = 
o(1)$. Similarly, one computes ${\mathcal L}_{X}g_{t} = {\mathcal L}_{X}g 
+ O(t^{n+1}) = 2t^{n}\widetilde \kappa + O(t^{n+1})$. Since $H \sim n$,
 using \eqref{e5.26} this gives
$$-\int_{S(t)}\langle {\mathcal L}_{X}(Hg_{t}), \widetilde \kappa \rangle dV_{S(t)} = 
-n\int_{S(t)}\langle {\mathcal L}_{X}g_{(n)}, \hat \kappa \rangle_{\gamma} dV_{\gamma} 
+ o(1).$$
Combining these computations and using \eqref{e5.26} again then gives \eqref{e5.27}. 

  Next, for the right side of \eqref{e5.25}, one has $A' = \frac{d}{ds}
(A_{g+s\widetilde \kappa}) = \frac{1}{2}({\mathcal L}_{N}\widetilde \kappa + 
{\mathcal L}_{N'}g) = \frac{1}{2}\nabla_{N}\widetilde \kappa + \widetilde \kappa 
+ O(t)$. Similarly, $(Hg_{t})' = H'g_{t} + H(g_{t})'$. The first term here, 
when paired with $(\delta^{*}X)^{T}$ and integrated, gives $O(t)$, 
while the second term is $n\widetilde \kappa$ to leading order. Hence 
\begin{equation}\label{e5.29}
2\int_{S(t)}\langle \tau'_{\widetilde \kappa}, (\delta^{*}X)^{T} \rangle dV_{g_{t}} = 
\int_{S(t)}\langle \nabla_{N}\widetilde \kappa - (2n-2)\widetilde \kappa, 
\kappa \rangle_{g_{t}} dV_{g_{t}} + O(t).
\end{equation}
A simple calculation shows that 
\begin{equation}\label{e5.30}
\int_{S(t)}\langle \nabla_{N}\widetilde \kappa - (2n-2)\widetilde \kappa, 
\kappa \rangle_{g_{t}} dV_{g_{t}} = 
\int_{S(t)}[{\tfrac{1}{2}}N(|\hat \kappa|^{2}) - (2n-2)|\hat \kappa|^{2}]
dV_{\gamma} + O(t),
\end{equation}
where the norms on the right are with respect to $\bar g$. The first term on 
the right in \eqref{e5.30} is $O(t)$, and using \eqref{e5.26} one obtains 
$$2\int_{S(t)}\langle \tau'_{\widetilde \kappa}, (\delta^{*}X)^{T}\rangle 
= -\frac{2n-2}{4}\int_{\partial M}|{\mathcal L}_{X}g_{(n)}|^{2}dV_{\gamma} 
+ o(1).$$
Comparing this with \eqref{e5.27} and using \eqref{e5.25}, it follows that on 
$\partial M$, 
$${\mathcal L}_{X}g_{(n)} = 0,$$
so that the flow of $X$ preserves both the boundary metric $\gamma$ and the 
$g_{(n)}$ term. The result then follows from the unique continuation result, 
[2, Corollary 4.4], analogous to Theorem 2.2.

{\endproof}

\section{Appendix: Corrections (October, 2012)}
\setcounter{equation}{0}

  In this Appendix, we point out and correct an error in the proof of Theorem 1.1 above. 
  
    The main problem with the proof of Theorem 1.1 is that the boundary conditions 
(5.9)-(5.10), i.e.
\be \label{1}
\int_{\dm}\<k(N), Y\> = 0,
\ee
\be \label{2}
\int_{\dm}(\d^{*}Y)(N, X) = 0,
\ee
cannot always be simultaneously enforced with the slice property (5.7), i.e.
\be \label{2a}
\int_{M}\<q, k \> \neq 0.
\ee
Namely, integration by parts in \eqref{2} gives
\be \label{3}
0 = \int_{\dm}(\d^{*}Y)(N, X) = \int_{\dm}(\d^{*}X)(N, Y) + \int_{M}\<\d \d^{*}Y, X\> 
- \int_{M}\<\d \d^{*}X, Y\>.
\ee
Suppose for example one chooses  $k = \d^{*}X \in K$, where $X$ is the given Killing field on $\dm$ 
preserving $H$. One has $\d k = \d \d^{*}X = 0$. It follows then from 
\eqref{1} and \eqref{3} with $k = \d^{*}X$ that 
$$0 = \int_{M}\< \d \d^{*}Y, X \> = -\int_{M}\<\d q, X\>,$$
where the second equality follows from the defining equation (5.8) for $Y$, i.e.~$\d \d^{*}Y = -\d q$. 
But
$$\int_{M}\<\d q, X\> = \int_{M}\<q, \d^{*}X\> + \int_{\dm}q(X,N) =  \int_{M}\<q, k\>,$$
since $q$ has compact support. Thus
\be \label{3a}
\int_{M}\<q, k\> = 0,
\ee
contradicting \eqref{2a}, i.e.~(5.7). This shows one cannot enforce both (5.9)-(5.10) given 
(5.7)-(5.8). 

\begin{remark}
{\rm  Referring to the discussion in the paragraphs following (5.11), suppose first $(M, g)$ 
has no Killing fields. Then the Dirichlet-to-Neumann map $R_{q}$ is surjective so that the boundary 
data $\d^{*}Y(N)$ for solutions of
\be \label{4}
\d \d^{*}Y = -\d q
\ee
can be arbitrarily prescribed on $\dm$. Let $\cS$ be the space of vector fields on $\dm$ which 
are $L^{2}$ orthogonal to $\{k(N)\}$, for all $k \in K$; $\cS$ is thus a subspace of codimension 
$k = dim K$ in the space $\chi(\dm)$ of all vector fields on $\dm$. The image $R_{q}(\cS)$ is 
then also of codimension $k$ in $\chi(\dm)$. Any $Y \in \cS$ such that $R_{q}(Y) \perp \<X\>$
then satisfies \eqref{1}-\eqref{2}. The space of such $Y$ is of codimension at most $k+1$, so 
obviously non-empty. Of course \eqref{2a} may also be satisfied, by choosing $q$ appropriately.  
It follows from from the discussion above in \eqref{3}-\eqref{3a} that in this situation
$$K \cap Im \d^{*} = 0,$$
i.e.~(5.22) holds. 

  On the other hand, if the space $\cK$ of Killing fields on $(M, g)$ is non-trivial, then the image 
$\cV_{q}$ of $R_{q}$ is a proper affine subspace of $\chi(\dm)$. When $R_{q}(\cV_{q})$ is an affine 
and not a linear subspace, it is then possible that there are no solutions of \eqref{4} satisfying 
\eqref{1} and \eqref{2}; this will be the case for instance if $X \in \cV_{0}$ but $X$ is orthogonal 
to $\cR_{0}(\cS)$. 
}
\end{remark}

   To resolve this problem, we proceed with essentially the same ideas and approach as before, 
but modify some of the details of the argument. 

   To begin,  consider $\w q$ now of the form
\be \label{5}
\w q = \psi g + D^{2}f,
\ee
where $\psi$ and $f$ are, for the moment, arbitrary smooth functions on $M$. Thus, we choose 
$q = \psi g$, no longer necessarily of compact support, and $Y = \nabla f$. Suppose 
that the Killing field $X$ is tangent to $\dm$, (as is the case for Theorem 1.1). We 
claim that \eqref{2} is automatically satisfied for this choice of $Y$. This follows from the 
following computation:
$$\int_{\dm}(\d^{*}Y)(N, X) = \int_{\dm}D^{2}f(N, X) = \int_{\dm}\<\nabla_{X}\nabla f, N\> 
= \int_{\dm}X(N(f)) - \<A(X), \nabla f\>$$ 
$$ = \int_{\dm}-div(X)N(f) - f(\d A)(X) + f\<A, \d^{*}X\> = 0.$$
The last equality follows since $div X = 0$ ($X$ is Killing on $\dm$), the divergence constraint 
$(\d A)(X) = -X(H) = 0$ and the Killing equation $\d^{*}X = 0$ on $\dm$. This proves 
the claim. 

   Moreover, referring to (5.19)-(5.21), since
$$\int_{\dm}q(N, X) = \int_{\dm}\psi \<X, N \> = 0,$$
(since $X$ is tangent to $\dm$) it follows that 
\be \label{5aa}
\int_{\dm}\w q(N, X) = 0,
\ee
i.e.~(5.21) holds.  

   The full argument of Section 5 (see in particular the analysis beginning with (5.14)) proceeds 
as before provided one has two properties. First the slice property, i.e.~the construction of a space 
$\w Q$ of smooth forms $\w q$ with $dim \w Q = dim K$ such that 
\be \label{5a}
Im \hat L \oplus \w Q = S_{2}^{m-2,\a}(M).
\ee
Since $Im \hat L = K^{\perp}$, cf.~the discussion preceding (5.6), this follows from the 
following property: for each $\w q \in \w Q$ there exists $k \in K$ such that 
\be \label{6}
\int_{M}\<\w q, k \> \neq 0,
\ee
so that no $\w q$ is orthogonal to $K$, cf.~(5.12). The second property is the divergence-free 
property (5.13), i.e.~$\d \w q = 0$, as in \eqref{4}. This requires 
$$-d\psi + \d(D^{2}f) = -d\psi - d\D f - Ric(df) = 0,$$
so that, since $Ric(df) = \l df$,
\be \label{7}
\D f = -\psi - \l f,
\ee
(up to an additive constant). In the following we choose $f$ to be arbitrary; \eqref{7} is then the 
defining equation for $\psi$ in \eqref{5}. It follows that for any $f$ and for $\psi$ satisfying 
\eqref{7}, one only needs to establish the slice property \eqref{6}. The condition \eqref{1} per 
se is dropped, since it was only used before to establish the slice property \eqref{6}. 

   Now computing \eqref{6} gives, since $\d k = 0$,
$$\int_{M}\< \w q, k\> = \int_{M}\psi tr k + \int_{\dm}\<k(N), df\> = 
\int_{M}\psi tr k + \int_{\dm}\d(k(N)^{T})f + k_{00}N(f),$$
where $k_{00} = k(N,N)$. We recall also that $k^{T} = \f \g$ on $\dm$. Set 
$$Z = \int_{\dm}\d(k(N)^{T})f + k_{00}N(f),$$
so 
\be \label{9}
\int_{M}\< \w q, k\> = \int_{M}\psi tr k + Z = -\int_{M}(\D f + \l f) tr k + Z.
\ee
On the other hand, since $\hat L(k) = 0$, cf.~(4.26), $tr k$ satisfies 
\be \label{9a}
\D tr k + \l tr k = 0.
\ee
Using this and integration by parts gives 
$$-\int_{M}\psi tr k = \int_{M}(\D f + \l f)tr k = \int_{\dm}N(f) tr k - N(tr k)f$$
$$ = \int_{\dm}N(f)k_{00} + N(f)n\f - N(tr k)f = R + Z + \int_{\dm}-\d(k(N)^{T})f + N(f)n\f - N(tr k)f,$$
$$ = Z + \int_{\dm}N(f)n \f - 2H'_{k}f - H\f f,$$
where we have used the standard formulas:
$$2H'_{k} = N(tr k) + 2\d(k(N)^{T}) - k_{00}H - N(k_{00}),$$
$$(\d k)(N) = -N(k_{00}) + \d(k(N)^{T}) + \<A, k\> - k_{00}H = 0.$$
Substituting this in \eqref{9} gives then the basic formula
\be \label{10}
\int_{M}\< \w q, k\> = -\int_{\dm}\f[nN(f) - Hf] - 2H'_{k}f.
\ee

   Now one requires only that \eqref{10} is non-zero, as in \eqref{6}. However, 
if \eqref{10} vanishes for all choices of $f$, then necessarily 
$$\f = 0 \ \ {\rm and} \ \  H'_{k} = 0.$$
Namely in \eqref{10} one can set $f = 0$ and $N(f)$ arbitrary on $\dm$ to obtain $\f = 0$; given 
this one can then set $N(f) = 0$ and $f$ arbitrary to obtain $H'_{k} = 0$. 

   Suppose for the moment that $\f = 0$, say for all $k \in K$. Then
\be \label{10a}
\int_{M}\<\w q, k\> = 2\int_{\dm}fH'_{k}.
\ee
For $k \in K$ one has of course $H'_{k} = 0$ and hence \eqref{10a} vanishes, for all choices 
of $\w q$, $k \in K$. Thus, the slice property \eqref{6} fails. 

  To circumvent this situation, consider the 
elliptic operator $\hat L$ as in (4.26) with boundary condition 
\be \label{11}
\hat L(h) = \ell, \ \d h = h_{0}, \ [h^{T}]_{0} = h_{1}, \ tr_{\s}A'_{h} - \<A, h \> = h_{2},
\ee
where $\s > 0$ is any smooth Riemannian metric on $\dm$. The first equation is on $M$, 
the last three on $\dm$. These boundary conditions are essentially exactly 
those considered in (3.10) with $B = \g$. Here $[h^{T}]_{0}$ denotes the usual 
tangential conformal equivalence class; the last term $\<A, h\>$ is lower order and is 
included only so that when $\s = \g$ one obtains, $tr_{\g}A'_{h} - \<A, h \> = H'_{h}$. 
Comparing \eqref{11} with (5.1) we are keeping the main factor $[h^{T}]_{0}$ the same, 
changing only the transverse scalar part from $H'_{h} = tr_{\g}A'_{h} - \<A, h\> = h_{2}$ to 
$tr_{\s}A'_{h} - \<A, h\> = h_{2}$. Note this scalar part played no role in Section 5. 

   By Proposition 3.2, \eqref{11} is an elliptic boundary value problem, 
of Fredholm index zero. The kernel $K_{\s}$ is finite dimensional and consists of 
forms $k_{\s}$ such that
$$\hat L(k_{\s}) = 0, \ \d k_{\s} = 0, \ [k_{\s}]_{0} = 0, \ tr_{\s}A'_{k_{\s}} - \<A, k_{\s}\> = 0.$$
Note that such $k_{\s}$ are still infinitesimal Einstein deformations. For general 
$\s$, the operator $\hat L$ on the space of symmetric forms $S_{0}^{m,\a}(M)$ satisfying vanishing 
boundary conditions in \eqref{11} may no longer have the property $K_{\s} \cap Im(\hat L) = 
0$, (cf.~\eqref{5a}). Thus, for simplicity, for the rest of the argument we assume 
that the metric $\s$ is (sufficiently) close to $\g$. In this case, $dim K_{\s}\leq 
dim K_{\g}$ and $K_{\s}$ is a slice to $Im \hat L$, i.e.
$$Im(\hat L) \oplus K_{\s} = S_{2}^{m-2,\a}(M).$$
In particular if $K_{\s} = 0$ for some $\s$, then the proof of Theorem 1.1 
carries over without change. Working exactly as before as in (5.14), it suffices then 
to find a ``good" slice $\w Q_{\s}$ for the kernel $K_{\s}$ in place of the original 
$K = K_{\g}$. In other words, \eqref{6} should hold, for $k = k_{\s} \in K_{\s}$. 

\medskip 

   Now it is straightforward to verify that all the computations above hold for 
$k = k_{\s} \in K_{\s}$ for any $\s$. In particular, \eqref{10} remains valid:
\be \label{12}
\int_{M}\< \w q, k_{\s}\> =  -\int_{\dm}\f_{k_{\s}}[nN(f) - Hf] - 2H'_{k_{\s}}f.
\ee

  For each choice of Riemannian metric $\s$ and each $k_{\s} \in K_{\s}$ one then has the boundary 
data $\f_{k_{\s}}$ and $H'_{k_{\s}}$. Whenever one can choose $\s$ and $k_{\s}$ such that 
$H'_{k_{\s}} \neq 0$, one can arrange that \eqref{12} is non-zero, by choosing $f$ and $N(f)$ 
suitably. Similarly, whenever $\f_{k_{\s}} \neq 0$ one can choose $f$ and $N(f)$ such that 
\eqref{12} is non-zero. 

 In general, for any $\s$ as above, consider the "reduced kernel" $\w K_{\s} \subset K_{\s}$ consisting 
of those $k_{\s} \in K_{\s}$ with $\f_{k_{\s}} = 0$. A basic property is that for $k_{\s} \in \w K_{\s}$,
$$H'_{k_{\s}} = 0 \Leftrightarrow k_{\s} \in \w K_{\g}.$$ 
Let $L_{\s}$ be a complement for $\w K_{\s}$ so that $K_{\s} = \w K_{\s} \oplus L_{\s}$. If $\ell_{j}$ 
is a basis for $L_{\s}$, then the boundary values $\f_{j}$ ($\ell_{j}^{T} = \f_{j}\g$ on $\dm$) are 
linearly independent. Hence, by choosing $f_{j} = 0$ and $N(f_{j})$ suitably, one obtains from \eqref{12} 
a space $\w Q_{L_{\s}}$ satisfying the slice property:
\be \label{13}
\int_{M}\< \w q, \ell\> \neq 0
\ee
for $\w q \in \w Q_{L_{\s}}$. 

   We now choose $\s$ as follows. For the ``original'' kernel $K = K_{\g}$ and its reduced kernel 
$\w K_{\g}$, choose $\s$ such that (as functions on $\dm$)
\be \label{14}
tr_{\s}A'_{k} \neq 0, 
\ee
for all $k \neq 0 \in \w K_{\g}$. (If $\w K_{\g} = 0$, so $K_{\g} = L_{\g}$, then \eqref{13} gives the required 
slice property for $K_{\s}$, $\s = \g$). If some $k_{\s} \in K_{\s}$ satisfies $k_{\s} = k\in \w K_{\g}$, then 
one has of course $tr_{\s}A'_{k} \neq 0$ by \eqref{14} but by definition of $K_{\s}$, $tr_{\s}A'_{k_{\s}} = 0$ 
(since $k_{\s} = 0$ on $\dm$) a contradiction. Thus $k_{\s} \notin \w K_{\g}$ for all $k_{\s}$, i.e.
$$K_{\s} \cap \w K_{\g} = 0,$$
for all $\s$ satisfying \eqref{14}. 

   Now the defining property of $\w K_{\g}$ (among forms satisfying \eqref{11} with $\ell = 0$, $h_{1} = 0$) 
is that $k \in \w K_{\g}$ if and only if $H'_{k} = 0$ and $\f = \f_{k} = 0$. Hence $k_{\s} \notin \w K_{\g}$ 
if and only if either $H'_{k_{\s}} \neq 0$ or $\f_{k_{\s}} \neq 0$. If the latter holds, 
then $k_{\s} \in L_{\s}$ and so \eqref{13} gives the slice property. If $\f_{k_{\s}} = 0$, then 
$k_{\s} \in \w K_{\s}$ but $H'_{k_{\s}} \neq 0$; if $k_{j}$ is a basis of $\w K_{\s}$, then the 
functions $H'_{k_{j}}$ are linearly independent. Thus again via \eqref{12} a suitable choice of $\{f_{j}\}$ 
gives the slice property \eqref{6} on $\w K_{\s}$ and so together with \eqref{13}, the slice property 
for all of $K_{\s}$. 

\medskip 

   To complete the proof, it thus suffices to prove there exists $\s > 0$ near $\g$ such that \eqref{14} 
holds. To do this, note first that for $k \in \w K_{\g}$, (so $\f_{k} = 0$ on $\dm$), $A'_{k} \neq 0$ on $\dm$. 
Namely, if $k^{T} = (A'_{k})^{T} = 0$ on $\dm$, it follows from the 
unique continuation theorem of [2] that $k = 0$ on $M$. Hence if $k_{j}$ is a basis for $\w K_{\g}$ 
then the symmetric forms $A'_{k_{j}}$ are linearly independent on $\dm$. 

    There are certainly many ways to prove the existence of $\s > 0$ for which \eqref{14} holds. 
One method is as follows. Note that \eqref{14} may be reformulated as: find a positive definite 
symmetric form $B$ such that
\be \label{14a}
tr_{\g}(BA'_{k}) \neq 0,
\ee
for all $k \neq 0 \in \w K_{\g}$.     

   Since each $A'_{k}$ is trace-free (since $tr A'_{k} = H'_{k} = 0$), each has a 
non-trivial positive part $(A'_{k})^{+}$ given by composing $A'_{k}$ with the projection onto the 
positive eigenspaces of $A'_{k}$. In particular, on any basis $k_{j}$ of $\w K_{\g}$, the forms 
$(A'_{k_{j}})^{+}$ are linearly independent; hence they are pointwise linearly independent on 
an open domain $\Omega \subset \dm$. To simplify the notation, set $A_{j}^{+} = (A'_{k_{j}})^{+}$ 
and $A_{j} = A'_{k_{j}}$. 

   Choose then points $p_{i} \in \Omega$, $1 \leq i \leq dim K$ with disjoint neighborhoods 
$U_{i} \subset \Omega$ and positive bump functions $\eta_{i}$ supported in $U_{i}$, with 
$\eta_{i}(p_{i}) = 1$. For the moment, set $B = \sum_{j}\eta_{j}A_{j}^{+}$, where for each 
$i$, the basis forms $\{A_{i}^{+}\}$ satisfy 
\be \label{15}
\<A_{i}^{+}, A_{j}\>(p_{i}) = 0, \ \ {\rm for \ all} \ \  j > i.
\ee
One constructs such a basis inductively as follows. At $p_{1}$ choose any basis $k_{i}$ of $\w K_{\g}$. 
Fix $k_{1}$ and $A_{1} = A'_{k_{1}}$ and then via the standard Gram-Schmidt process, construct
the basis forms $k_{j}$, $j \geq 2$ satisfying \eqref{15} at $p_{1}$. Next in the space spanned 
by $\{k_{j}\}$, $j \geq 2$, repeat the process at $p_{2}$, starting with $A_{2}$ and constructing 
forms $k_{j}$, $j \geq 3$ satisfying \eqref{15} at $p_{2}$. One continues inductively in this 
way through to the last point. Note that a different basis of $K$ is thus used at each point $p_{i}$. 
At any given $p_{r}$ one has
\be \label{15a}
tr_{\g}(BA_{k})(p_{r}) = \<B, A_{k}\>(p_{r})= \sum_{i,j}\eta_{i}c_{j}\<A_{i}^{+}, A_{j}\>(p_{r}),
\ee
where $k = \sum c_{j}k_{j}$ in the basis associated to $p_{r}$. 

   Now suppose \eqref{14a} fails for some $k$, so that $tr_{\g}(BA'_{k}) = 0$ on $\dm$. 
Evaluating \eqref{15a} at $p_{1}$ gives, by \eqref{15},
$$tr_{\g}(BA'_{k})(p_{1}) = c_{1}|A_{1}^{+}|^{2} = 0,$$
so that $c_{1} = 0$. From this, and from the construction of the basis at $p_{2}$, one has 
similarly
$$tr_{\g}(BA'_{k})(p_{2}) = c_{2}|A_{2}^{+}|^{2} = 0,$$
so that $c_{2} = 0$. Continuing in this way, it follows that $c_{r} = 0$ for all $r$, and hence 
(again by the construction of the bases) $k = 0$. This establishes the property \eqref{14a} 
for $B$ as above.

   Finally, note that for $B' = Id$, $tr_{\g}(B'A'_{k}) = tr_{\g}(A'_{k}) = H'_{k} = 0$, for all 
$k \in \w K_{\g}$. Also, for $B$ as above, on the unit sphere in $\w K_{\g}$ the space of functions 
$tr_{\g}(BA'_{k})$ is compact, and so bounded away from the zero function. Hence, choosing $\e$ 
sufficiently small and replacing $B$ by $Id + \e B$ gives a smooth metric $\s > 0$, close to $\g$ 
on $\dm$, satisfying \eqref{14}. 

     Given the changes above, the rest of the proof of Theorem 1.1 remains the same. 
Briefly, the construction above gives the existence of a slice $\w Q_{\s}$ as in \eqref{5a} 
consisting of divergence-free forms $\w q$ which satisfy \eqref{5aa}. One constructs then 
the operator $\w L = \hat L + \pi_{\w Q_{\s}}$ as in (5.14) and proceeds exactly as 
before to complete the proof. 
     
{\endproof}     

\bigskip 

  Essentially the same method holds for the proof of Proposition 5.2 generalizing Theorem 
1.1. to the case $X$ is an infinitesimal isometry at $\dm$ preserving the mean curvature $H$, 
i.e.~  $(\d^{*}X)^{T} = 0$, $H'_{\d^{*}X} = 0$, with $X$ not necessarily tangent to $\dm$.    
  
   Thus, define $\w q$ as in \eqref{5} with $\psi$ given by \eqref{7}, with $f$ free. Then as before 
$\d \w q = 0$ and as in (5.23)-(5.24) one needs the relation 
\be \label{21}
\int_{\dm}\hat E'_{h}(N,X) = 0,
\ee
to hold. By construction, cf.~(5.19), $E'_{h} = -\w q = -\psi g - D^{2}f$, so that \eqref{21} is 
equivalent to 
\be \label{22}
\int_{\dm}\<N, X\>(\D f + \l f) - D^{2}f(N,X) = 0.
\ee
To compute the second term in \eqref{22}, one has 
$$\int_{\dm}\<\nabla_{X}\nabla f, N\> = \int_{\dm}NN(f)\<X, N\> + X^{T}\<\nabla f, N\> - 
\<\nabla f, \nabla_{X^{T}}N\>,$$
$$= \int_{\dm}NN(f)\<X, N\> - div(X^{T})N(f) - A(X^{T}, \nabla f) = \int_{\dm}NN(f)\<X, N\> - 
div(X^{T})N(f) - f\d(A(X^{T})) $$
$$= \int_{\dm}NN(f)\<X, N\> - div(X^{T})N(f) + f\<A, \d^{*}X^{T}\> + fdH(X^{T}),$$
where we have used the fact that $\d A (X^{T}) = -dH(X^{T})$. Since $(\d^{*}X)^{T} = 0$, one 
has $\d^{*}X^{T} + \<X, N\>A = 0$, so that $div(X^{T}) = - \<X, N\>H$. It follows that
\be \label{23}
\int_{\dm}D^{2}f(N, X) = \int_{\dm}\<X, N\>[NN(f) + HN(f) - f|A|^{2}] + fX^{T}(H).
\ee

   On the other hand, for the first term in \eqref{22} one has $\D f = \D_{\dm}f + HN(f) + NN(f)$ so that 
setting $\nu = \<X, N\>$, 
\be \label{24}
\int_{\dm}\nu(\D f + \l f) = 
\int_{\dm}f \D_{\dm}\nu + HN(f)\nu + NN(f)\nu + \l f \nu.
\ee
Subtracting \eqref{23} from \eqref{24} gives
$$\int_{\dm}\hat E'_{h}(N,X) = \int_{\dm}f[\D \nu + (|A|^{2} + \l)\nu - X^{T}(H)] = -\int_{\dm}fH'_{X} = 0,$$
since $X$ preserves the mean curvature; the second equality here is exactly the formula for the variation 
of the mean curvature in the direction $X$. 
   
   This completes the corrected proof of Proposition 5.2. 

{\endproof}

   Again, with the modifications above, the proof of Theorem 5.4 proceeds just as before. 

{\endproof}

\bibliographystyle{plain}

\begin{thebibliography}{WWW}

\footnotesize


\bibitem [1]{AG} M. Akbar and G. Gibbons, Ricci-flat metrics with $U(1)$ action and the 
Dirichlet boundary-value problem in Riemannian quantum gravity and isoperimetric 
inequalities, Class. Quantum Gravity, {\bf 20}, (2003), 1787-1822. 

\bibitem [2]{AH} M. Anderson and M. Herzlich, Unique continuation results for Ricci 
curvature and applications, Jour. Geometry \& Physics, {\bf 58}, (2008), 179-207, 
arXiv:0710.1305 [math.DG]; see also Erratum, {\it ibid}, {\bf 60}, (2010), 1062-1067.  

\bibitem [3]{An1} M. Anderson, On boundary value problems for Einstein metrics, 
Geometry \& Topology, {\bf 12}, (2008), 2009-2045, arXiv: math.DG/0612647. 

\bibitem [4]{AK} M. Anderson and M. Khuri, The static extension problem in 
general relativity, (preprint, Sept, 2009), arXiv:0909.4550 (math.DG). 

\bibitem [5]{Be} A. Besse, Einstein Manifolds, Springer Verlag, New York, (1987).

\bibitem [6]{Bi} O. Biquard, Continuation unique a partir de l'infinie conforme pour les 
metriques d'Einstein, Math. Research Lett., {\bf 15}, (2008), 1091-1099. 

\bibitem [7]{Bo} A. Borisenko, Isometric immersions of space forms into Riemannian 
and pseudo-Riemannian spaces of constant curvature, Russian Math. Surveys, {\bf 56}, 
(2001), 425-497. 

\bibitem [8]{GT} D. Gilbarg and N. Trudinger, Elliptic Partial Differential Equations 
of Second Order, Springer Verlag, New York, (1983). 

\bibitem [9]{Ha} S. Hawking, Euclidean quantum gravity, in: Recent Developments in 
Gravitation, Cargese Lectures, eds. M. Levy and S. Deser, (Plenum 1978), 145-173. 

\bibitem [10]{Ki} Y. Kitagawa, Deformable flat tori in $S^{3}$ with constant mean 
curvature, Osaka Math. Jour., {\bf 40}, (2003), 103-119. 

\bibitem [11]{KN} S. Kobayashi and K. Nomizu, Foundations of Differential Geometry, 
vol. 1, Wiley-Interscience, New York, (1963). 

\bibitem [12]{M} C. Morrey, Jr. Multiple Integrals in the Calculus of Variations, 
Grundlehren Series 130, Springer Verlag, New York, (1966). 

\bibitem [13]{P} P. Petersen, Riemannian Geometry, Graduate Texts in Mathematics, 
vol. 171, Springer Verlag, New York, (1997). 

\bibitem [14]{Pi} U. Pinkall, Hopf tori in $S^{3}$, Inventiones Math., {\bf 81}, 
(1985), 379-386. 

\bibitem [15]{R} H. Rosenberg, (private communication). 

\bibitem [16]{Sch} J.M. Schlenker, Einstein manifolds with convex boundaries, 
Comm. Math. Helv., {\bf 76}, (2001), 1-28. 

\bibitem [17]{Sh} J.P. Sha, $p$-convex Riemannian manifolds, Inventiones Math., 
{\bf 83}, (1986), 437-447. 

\bibitem [18]{Sp} M. Spivak, A Comprehensive Introduction to Differential Geometry, 
Vol. III, V, $2^{\rm nd}$ Edition, Publish or Perish, Inc., Berkeley, (1979). 

\bibitem [19]{Yo} J. W. York, Jr., Black-hole thermodynamics and the Euclidean 
Einstein action, Phys. Rev. D, {\bf 33}, (1986), 2092-2099. 



\end{thebibliography}

\smallskip
\noindent
\address{Department of Mathematics\\
Stony Brook University\\
Stony Brook, N.Y. 11794-3651}

\noindent
E-mail: anderson@math.sunysb.edu

\end{document}